\documentclass[11pt]{article}
\usepackage{newlfont,amsfonts,amssymb,
amsmath,amsthm,amsgen,amscd}
\def\a{\alpha}
\def\b{\beta}

\def\d{\delta}
\def\e{\varepsilon}

\def\g{\gamma}
\def\h{\eta}
\def\i{\iota}

\def\l{\lambda}

\def\t{\tau}
\def\u{\upsilon}

\def\G{\Gamma}

\def\O{\Omega}

\chardef\tempcat=\the\catcode`\@
\catcode`\@=11
\def\cyracc{\def\u##1{\if \i##1\accent"24 i
    \else \accent"24 ##1\fi }}
\newfam\cyrfam

\DeclareFontFamily{OT1}{msb}{}{}
\DeclareFontShape{OT1}{msb}{m}{n}
 {  <5> <6> <7> <8> <9> <10> gen * msbm
      <10.95><12><14.4><17.28><20.74><24.88>msbm10}{}
\DeclareMathAlphabet{\bubble}{OT1}{msb}{m}{n}

\def\bR{{\mathbb R}}

\newfont{\goth}{eufm10 scaled \magstep1}

\newfont{\mcal}{eusm10 scaled \magstep1}

\def\ch{\mbox{\mcal H}}

\def\p{\partial}

\def\square{\kern1pt\vbox
            {\hrule height 0.6pt\hbox{\vrule width 0.6pt\hskip 3pt
 \vbox{\vskip 6pt}\hskip 3pt\vrule width 0.6pt}\hrule height 0.6pt}\kern1pt}

\def\hath{\widehat{\mathfrak{h}}}

\def\t{\tilde}

\newtheorem{Th}{Theorem}
\newtheorem{Prop}{Proposition}
\newtheorem{Cor}{Corollary}
\newtheorem{Lem}{Lemma}
\newtheorem{Def}{Definition}

\def\bt{\begin{Th}}
\def\et{\end{Th}}
\def\bp{\begin{Prop}}
\def\ep{\end{Prop}}
\def\bc{\begin{Cor}}
\def\ec{\end{Cor}}
\def\bl{\begin{Lem}}
\def\el{\end{Lem}}
\def\bd{\begin{Def}}
\def\ed{\end{Def}}

\def\n{\nabla}
 \def\ot{\otimes}

\def\arr{\begin{array}}
\def\ea{\end{array}}
\def\bea{\begin{eqnarray}}
\def\eea{\end{eqnarray}}
\def\bean{\begin{eqnarray*}}
\def\eean{\end{eqnarray*}}

\catcode`@=11
\@addtoreset{equation}{section}
\catcode`@=12

\def\arctg{\mathrm{arctg\,}}
\def\arcth{\mathrm{artanh\,}}
\def\sh{\mathrm{sinh\,}}
\def\ch{\mathrm{cosh\,}}
\def\tg{\mathrm{tan\,}}
\def\th{\mathrm{tanh\,}}

\def\h{\hat}
\def\hol{\mathfrak{hol}}
\def\so{\mathfrak{so}}

\begin{document}
\newcommand{\N}{\ensuremath{\mathbb{N}}}
\newcommand{\Z}{\ensuremath{\mathbb{Z}}}
\newcommand{\C}{\ensuremath{\mathbb{C}}}
\newcommand{\K}{\ensuremath{\mathbb{K}}}
\renewcommand{\O}{\ensuremath{\mathcal{O}}}
\newcommand{\R}{\ensuremath{\mathbb{R}}}

\newcommand{\bcase}{\begin{case}}
\newcommand{\ecase}{\end{case}}
\newcommand{\setcase}{\setcounter{case}{0}}
\newcommand{\bclaim}{\begin{claim}}
\newcommand{\eclaim}{\end{claim}}
\newcommand{\setclaim}{\setcounter{claim}{0}}
\newcommand{\bstep}{\begin{step}}
\newcommand{\estep}{\end{step}}
\newcommand{\setstep}{\setcounter{step}{0}}
\newcommand{\bhlem}{\begin{hlem}}
\newcommand{\ehlem}{\end{hlem}}
\newcommand{\sethlem}{\setcounter{hlem}{0}}

\newcommand{\bleer}{\begin{leer}}
\newcommand{\eleer}{\end{leer}}
\newcommand{\bde}{\begin{de}}
\newcommand{\ede}{\end{de}}
\newcommand{\ul}{\underline}
\newcommand{\ol}{\overline}
\newcommand{\tbf}{\textbf}
\newcommand{\mc}{\mathcal}
\newcommand{\mb}{\mathbb}
\newcommand{\mf}{\mathfrak}
\newcommand{\bs}{\begin{satz}}
\newcommand{\es}{\end{satz}}
\newcommand{\btheo}{\begin{theo}}
\newcommand{\etheo}{\end{theo}}
\newcommand{\bfolg}{\begin{folg}}
\newcommand{\efolg}{\end{folg}}
\newcommand{\blem}{\begin{lem}}
\newcommand{\elem}{\end{lem}}
\newcommand{\bnote}{\begin{note}}
\newcommand{\enote}{\end{note}}
\newcommand{\bprf}{\begin{proof}}
\newcommand{\eprf}{\end{proof}}
\newcommand{\be}{\begin{eqnarray*}}
\newcommand{\ee}{\end{eqnarray*}}
\newcommand{\eeqa}{\end{eqnarray}}
\newcommand{\beqa}{\begin{eqnarray}}
\newcommand{\bi}{\begin{itemize}}
\newcommand{\ei}{\end{itemize}}
\newcommand{\bnum}{\begin{enumerate}}
\newcommand{\enum}{\end{enumerate}}
\newcommand{\la}{\langle}
\newcommand{\ra}{\rangle}
\newcommand{\eps}{\epsilon}
\newcommand{\ve}{\varepsilon}
\newcommand{\vp}{\varphi}
\newcommand{\lra}{\longrightarrow}
\newcommand{\Lra}{\Leftrightarrow}
\newcommand{\Ra}{\Rightarrow}
\newcommand{\sub}{\subset}
\newcommand{\ems}{\emptyset}
\newcommand{\sms}{\setminus}
\newcommand{\ints}{\int\limits}
\newcommand{\sums}{\sum\limits}
\newcommand{\lims}{\lim\limits}
\newcommand{\bcup}{\bigcup\limits}
\newcommand{\bcap}{\bigcap\limits}
\newcommand{\beq}{\begin{equation}}
\newcommand{\eeq}{\end{equation}}
\newcommand{\einhalb}{\frac{1}{2}}
\newcommand{\rr}{\mathbb{R}}
\newcommand{\rn}{\mathbb{R}^n}
\newcommand{\ccc}{\mathbb{C}}
\newcommand{\cn}{\mathbb{C}^n}
\newcommand{\M}{{\cal M}}
\newcommand{\drehgleich}{\mbox{\begin{rotate}{90}$=$  \end{rotate}}}
\newcommand{\turngleich}{\mbox{\begin{turn}{90}$=$  \end{turn}}}
\newcommand{\turnsimeq}{\mbox{\begin{turn}{270}$\simeq$  \end{turn}}}
\newcommand{\vf}{\varphi}
\newcommand{\earr}{\end{array}\]}
\newcommand{\barr}{\[\begin{array}}
\newcommand{\bvec}{\left(\begin{array}{c}}
\newcommand{\evec}{\end{array}\right)}
\newcommand{\sumk}{\sum_{k=1}^n}
\newcommand{\sumi}{\sum_{i=1}^n}
\newcommand{\suml}{\sum_{l=1}^n}
\newcommand{\sumj}{\sum_{j=1}^n}
\newcommand{\sumij}{\sum_{i,j=1}^n}
\newcommand{\suminf}{\sum_{k=0}^\infty}
\newcommand{\inv}{\frac{1}}
\newcommand{\wzbw}{\hfill \eprf\\[0.2cm]}
\newcommand{\lag}{\mathfrak{g}}
\newcommand{\lan}{\mathfrak{n}}
\newcommand{\lah}{\mathfrak{h}}
\newcommand{\laz}{\mathfrak{z}}
\newcommand{\+}{\oplus}
\newcommand{\lx}{\ltimes}
\newcommand{\rrn}{\mathbb{R}^n}
\newcommand{\laso}{\mathfrak{so}}
\newcommand{\lason}{\mathfrak{so}(n)}
\newcommand{\lagl}{\mathfrak{gl}}
\newcommand{\lasl}{\mathfrak{sl}}
\newcommand{\lasp}{\mathfrak{sp}}
\newcommand{\lasu}{\mathfrak{su}}
\newcommand{\w}{\omega}
\newcommand{\pmh}{{\cal P}(M,h)}
\newcommand{\s}{\sigma}
\newcommand{\deri}{\frac{\partial}}
\newcommand{\ddxi}{\frac{\partial}{\partial x_i}}
\newcommand{\ddxj}{\frac{\partial}{\partial x_j}}
\newcommand{\ddxk}{\frac{\partial}{\partial x_k}}
\newcommand{\ddxl}{\frac{\partial}{\partial x_l}}
\newcommand{\ddxp}{\frac{\partial}{\partial x_p}}
\newcommand{\ddxq}{\frac{\partial}{\partial x_q}}
\newcommand{\ddxr}{\frac{\partial}{\partial x_r}}
\newcommand{\ddxin}{\frac{\partial}{\partial x_{i+n}}}
\newcommand{\ddxjn}{\frac{\partial}{\partial x_{j+n}}}
\newcommand{\ddxkn}{\frac{\partial}{\partial x_{k+n}}}
\newcommand{\ddxln}{\frac{\partial}{\partial x_{l+n}}}
\newcommand{\ddxpn}{\frac{\partial}{\partial x_{p+n}}}
\newcommand{\ddxqn}{\frac{\partial}{\partial x_{q+n}}}
\newcommand{\hnab}{\hut{\nabla}}
\newcommand{\ddt}{\frac{\partial}{\partial t}}
\newcommand{\ddr}{\frac{\partial}{\partial r}}
\newcommand{\ddi}{\frac{\partial}{\partial y_i}}
\newcommand{\ddj}{\frac{\partial}{\partial y_j}}
\newcommand{\ddk}{\frac{\partial}{\partial y_k}}
\newcommand{\ddp}{\frac{\partial}{\partial p_i}}
\newcommand{\ddq}{\frac{\partial}{\partial q_i}}
\newcommand{\xz}{^{(x,z)}}
\newcommand{\mh}{(M,h)}
\newcommand{\wxz}{W_{(x,z)}}
\newcommand{\qmh}{{\cal Q}(M,h)}
\newcommand{\bbem}{\begin{bem}}
\newcommand{\ebem}{\end{bem}}
\newcommand{\bbez}{\begin{bez}}
\newcommand{\ebez}{\end{bez}}
\newcommand{\bbsp}{\begin{bsp}}
\newcommand{\ebsp}{\end{bsp}}
\newcommand{\pr}{pr_{\lason}}
\newcommand{\hut}{\widehat}
\newcommand{\huts}{\widehat{\s}}
\newcommand{\whut}{\w^{\huts}}
\newcommand{\bhg}{{\cal B}_H(\lag)}
\newcommand{\aaa}{\alpha}
\newcommand{\bb}{\beta}
\newcommand{\laa}{\mf{a}}
\newcommand{\lam}{\lambda}
\newcommand{\LL}{\Lambda}
\newcommand{\Š}{\alpha}
\newcommand{\W}{\Omega}
\newcommand{\esel}{\ensuremath{\mathfrak{sl}(2,\ccc)}}
\newcommand{\kg}{{\cal K}(\lag)}
\newcommand{\bg}{{\cal B}_h(\lag)}
\newcommand{\kk}{  \mathbb{K}}
\newcommand{\xy}{[x,y]}
\newcommand{\perdef}{$\stackrel{\text{\tiny def}}{\iff}$}
\newcommand{\eqdef}{\stackrel{\text{\tiny def}}{=}}
\newcommand{\lai}{\mf{i}}
\newcommand{\lar}{\mf{r}}
\newcommand{\Dim}{\mathsf{dim\ }}
\newcommand{\im}{\mathsf{im\ }}
\newcommand{\Ker}{\mathsf{ker\ }}
\newcommand{\trace}{\mathsf{trace\ }}
\newcommand{\grad}{\mathsf{grad}}
\newcommand{\wt}{\widetilde}
\newcommand{\tnab}{\widetilde{\nabla}}
\newcommand{\tem}{\widetilde{M}}
\newcommand{\nabt}{\nabla^{\cal T}}
\newcommand{\ro}{\mathsf{P}}
\newcommand{\lecturecount}{\begin{center}{\sf  (Lecture \Roman{lecturenr})}\end{center}\addtocontents{toc}{{\sf  (Lecture \Roman{lecturenr})}} \refstepcounter{lecturenr}
}
\newcommand{\T}{{\cal T}}
\newcommand{\cur}{{\cal R}}
\newcommand{\pd}{{\cal P}}

\newcommand{\inter}{\makebox[11pt]{\rule{6pt}{.3pt}\rule{.3pt}{5pt}}}


\theoremstyle{definition}
\newtheorem{de}{Definition}[section]
\newtheorem{bem}{Remark}[section]
\newtheorem{bez}{Notation}[section]
\newtheorem{bsp}{Example}[section]
\theoremstyle{plain}
\newtheorem{lem}{Lemma}[section]
\newtheorem{satz}{Proposition}[section]
\newtheorem{folg}{Corollary}[section]
\newtheorem{theo}{Theorem}[section]

\def\formtmaintwo{Let $(M,g)$ be a  pseudo-Riemannian manifold.
If the holonomy group of the metric cone over $(M,g)$  admits a non-degenerate
invariant subspace,
then $\hol(M,g)=\so(p,q)$, where $(p,q)$ is the signature of the metric $g$.}

\def\formtmainfive{
Let $(M,g)$ be a  pseudo-Riemannian manifold  and $(\widehat{M}=\bR^+\times M,\widehat{g}=dr^2+r^2g)$  the cone over $M$.
Suppose that the holonomy group of $(\widehat{M},\widehat{g})$ admits a non-degenerate proper invariant  subspace.
Then  there exists a dense open submanifold $U_1\subset M$ such that each point $x\in U_1$ has an open neighborhood $W\subset U_1$
that satisfies one of the following conditions
\begin{itemize}
\item[(1.)] For $W$ we have a decomposition $$W=(a,b)\times N_1\times N_2,\quad (a,b)\subset \left(0,\frac{\pi}{2}\right)$$
and for the metric $g|_{W}$ we have
$$g|_W=ds^2+ \cos^2(s) g_1+\sin^2(s)g_2,$$
where $(N_1,g_1)$ and $(N_2,g_2)$ are pseudo-Riemannian manifolds;

Moreover, any point  $(r,x)\in\bR^+\times W\subset\widehat{M}$  has a neighborhood of the form
$$((a_1,b_1)\times N_1)\times ((a_2,b_2)\times N_2),\quad (a_1,b_1),(a_2,b_2)\subset\bR^+$$
with the metric $$(dt_1^2+t_1^2g_1)+(dt_2^2+t_2^2g_2).$$

\item[(2.)]  For $W$ we have a decomposition $$W=(a,b)\times N_1\times N_2,\quad (a,b)\subset\bR^+$$
and for the metric $g|_{W}$ we have
$$g|_W=-ds^2+ \ch^2(s) g_1+\sh^2(s)g_2,$$
where $(N_1,g_1)$ and $(N_2,g_2)$ are pseudo-Riemannian manifolds.

Moreover, any point  $(r,x)\in\bR^+\times W\subset\widehat{M}$  has a neighborhood of the form
$$((a_1,b_1)\times N_1)\times ((a_2,b_2)\times N_2),\quad (a_1,b_1),(a_2,b_2)\subset\bR^+$$
with the metric $$(dt_1^2+t_1^2g_1)+(-dt_2^2+t_2^2g_2).$$
\end{itemize}}

\def\formtmainfore{
Let $(M,g)$ be a complete  pseudo-Riemannian  manifold of dimension $\ge 2$ with decomposable holonomy $\hath$ of the cone
$\widehat{M}$. Then there exists an open dense submanifold $M'\subset M$ such that each connected component of $M'$   is
isometric to a pseudo-Riemannian manifold of the form
\begin{itemize}
\item[(1)] a pseudo-Riemannian manifold $M_1$ of constant sectional curvature 1 or
\item[(2)] a pseudo-Riemannian manifold $M_2=\bR^+\times N_1\times N_2$
with the metric
$$-ds^2+ \ch^2(s) g_1+\sh^2(s)g_2,$$
where $(N_1,g_1)$ and $(N_2,g_2)$ are pseudo-Riemannian manifolds
and $(N_2,g_2)$ has constant sectional curvature $-1$ or $\dim N_2 \le 1$.

Moreover, the cone $\hut M_2$ is isometric to the open subset $\{ r_1>r_2\}$
in   the product of the space-like cone $(\bR^+\times N_1, dr^2 + r^2 g_1)$
over $(N_1,g_1)$ and the time-like cone
$(\bR^+ \times N_2,-dr^2 + r^2 g_2)$ over $(N_2,g_2)$.
\end{itemize}
}

\def\formtmainthree{Let $(M,g)$ be a compact and complete pseudo-Riemannian
manifold of dimension $\ge 2$
with decomposable holonomy group $\widehat{H}$ of the
cone $\widehat{M}$. Then $(M,g)$ has constant curvature 1 and the cone is flat.
}

\def\formtmainsix3Sasaki{
Let $(M,g)$ be a pseudo-Riemannian manifold. There
is a one-to-one correspondence between para-3-Sasakian structures
$(M,g,T_1,T_2,T_3)$ on $(M,g)$ and para-hyper-K\"ahler structures
$(\widehat{M},\widehat{g},\widehat{J}_1,\widehat{J}_2,\widehat{J}_3=\widehat{J}_1\widehat{J}_2)$
on the cone
$(\widehat{M},\widehat{g})$. The correspondence is given by
$T_\alpha \mapsto \widehat{J}_\alpha := \widehat{\nabla}T_\alpha$.
}

\def\formtmainseven{\begin{itemize}\item[1.]
Let $(M,g)$ be a Lorentzian manifold of signature $(+,\cdots,+,-)$ or a negative definite Riemannian manifold
and $(\widehat{M} = \bR^+ \times M, \widehat{g})$ the cone over $M$ equipped
with the Lorentzian metric $\widehat{g} = dr^2 + r^2 g$ (of signature $(+,\cdots,+,-)$ or $(+,-,\cdots,-)$).
Suppose that the cone $(\widehat{M},\widehat{g})$ admits a parallel distribution of isotropic lines.
If $M$ is simply connected, then $(\widehat{M},\widehat{g})$ admits a non-zero  parallel light-like vector field.

\item[2.] Let $(M,g)$ be a negative definite Riemannian manifold
and $(\widehat{M}, \widehat{g})$ the cone over $M$ equipped with the Lorentzian metric of signature  $(+,-,\cdots,-)$.
Suppose that the cone $(\widehat{M},\widehat{g})$ admits a  non-zero  parallel light-like vector field, then each point
$x\in M$ has a neighbourhood of the form $$M_0=(a,b)\times
N,\,\,\,\,\,a\in\bR\cup\{-\infty\},\,\,b\in\bR\cup\{+\infty\},\,\,a<b,$$ and for the metric $g|_{M_0}$ we have
$$g|_{M_0}=-ds^2+e^{-2s}g_N,$$ where $(N,g_N)$ is a negative definite Riemannian manifold.
If the holonomy algebra $\hol(\widehat{M_0},\widehat{g|_{M_0}})$ of the manifold $(\widehat{M_0},\widehat{g|_{M_0}})$ is
indecomposable, then $$\hol(\widehat{M_0},\widehat{g|_{M_0}})\cong\hol(N,g_N)\ltimes\bR^{\dim N}.$$ If the manifold
$(M,g)$ is complete, then $M_0=M$, $(a,b)=\bR$ and $(N,g_N)$ is complete.

\item[3.] Let $(M,g)$ be a Lorentzian  manifold
and $(\widehat{M},\widehat{g})$ the cone over $M$ equipped with the Lorentzian metric. Suppose that the cone
$(\widehat{M},\widehat{g})$ admits a  non-zero  parallel light-like vector field, then there exist disjoint open subspaces
$\{W_i\}_{i\in I}\subset M$ such that the open subspace $\cup_{i\in I} W_i\subset M$ is dense. Any point $x$ of each $W_i$
has an open neighbourhood of the form
 $$U_i=(a,b)\times N\subset W_i,\,\,\,\,\,a\in\bR\cup\{-\infty\},\,\,b\in\bR\cup\{+\infty\},\,\,a<b,$$ and for the metric
 $g|_{U_i}$ we have
$$g|_{U_i}=-ds^2+e^{-2s}g_N,$$ where $(N,g_N)$ is a Riemannian manifold.
If the holonomy algebra $\hol(\widehat{U_i},\widehat{g|_{U_i}})$ of the manifold $(\widehat{U_i},\widehat{g|_{U_i}})$ is
indecomposable, then  $$\hol(\widehat{U_i},\widehat{g|_{U_i}})\cong\hol(N,g_N)\ltimes\bR^{\dim N}.$$ If  the manifold
$(M,g)$ is complete, then $(a,b)=\bR$ and $U_i=W_i$.
\end{itemize}
}

\bibliographystyle{alpha}


\title{Cones over pseudo-Riemannian manifolds and their holonomy}

\author{D.\,V.\,Alekseevsky\thanks{{\sc
The University of Edinburgh and Maxwell Institute for Mathematical Sciences, JCMB, The Kings buildings, Edinburgh, EH9
3JZ, UK}, {\tt D.Aleksee@ed.ac.uk}},\hspace{.1cm} V.\,Cort\'es\thanks{{\sc  Department Mathematik und Zentrum f\"ur
Mathematische Physik, Universit\"at Hamburg, Bundesstra{\ss}e 55, D-20146 Hamburg, Germany}, {\tt
cortes@math.uni-hamburg.de}, {\tt leistner@math.uni-hamburg.de}},\hspace{.1cm} A.\,S.\,Galaev\thanks{{\sc  Masaryk University, Faculty of Science, Department
of Mathematics, Jan\'a\v ckovo n\' am. 2a, 60200 Brno, Czech Republic}, {\tt galaev@math.muni.cz}
\newline {\em Date:} July 25, 2007},\hspace{.1cm}
T.\,Leistner${}^\dag$
}
\date{ }
\maketitle
 \begin{abstract}
By a classical theorem of Gallot (1979), a Riemannian cone over a
complete Riemannian manifold is either flat or has irreducible holonomy.
We consider metric cones with reducible holonomy
over pseudo-Riemannian manifolds. First we describe
the local structure of the base of the cone
when the holonomy of the cone is decomposable.
For instance, we find that the holonomy algebra of the base is always
the full pseudo-orthogonal Lie algebra. One of the global results
is  that a cone over a
compact and complete pseudo-Riemannian manifold is either flat or has
indecomposable holonomy. Then we analyse the case when the cone has
indecomposable but reducible holonomy, which means that it admits  a
parallel isotropic distribution. This analysis is carried out, first  in the case where the
cone admits two complementary distributions and, second for
Lorentzian cones. We show that
the first case occurs precisely when the local geometry of the
base manifold is para-Sasakian and that of the cone is para-K\"ahlerian.
For Lorentzian cones we
get a complete description of the possible (local) holonomy
algebras in terms of the metric of the base manifold.\\[2mm]
\end{abstract}

\tableofcontents



\section{Introduction}
Let $(M,g)$ be a (connected) pseudo-Riemannian manifold of signature
$$(p,q)=(-,\cdots ,-,+,\cdots ,+).$$ We denote by $H$ the
holonomy group of $(M,g)$ and by $\mathfrak{h} \subset
\mathfrak{so}(V)$ ($V=T_pM$, $p\in M$)
the corresponding holonomy algebra.

We say that $\mathfrak{h}$ is \emph{decomposable} if
$V$ contains a proper non-degenerate $\mathfrak{h}$-invariant
subspace. By Wu's theorem \cite{wu64} this means that $M$ is
locally decomposed as a  product of two pseudo-Riemannian manifolds.
In the opposite case $\mathfrak{h}$ is called \emph{indecomposable}.
We say that $\mathfrak{h}$ is \emph{reducible} if it preserves
a (possibly degenerate) proper subspace of $V$.

The holonomy algebra $\mathfrak{h}$ is of exactly one
of the following types:
\begin{itemize}
\item[(i)] decomposable,
\item[(ii)] reducible indecomposable or
\item[(iii)] irreducible.
\end{itemize}

Let $(\widehat{M}=\bR^+\times M,\widehat{g}=dr^2+r^2g)$ be the
(space-like) metric cone over $(M,g)$. We denote by $\widehat{H}$ the
holonomy group of $(\widehat{M},\widehat{g})$ and by $\widehat{\mathfrak{h}}
\subset \mathfrak{so}(\widehat{V})$ ($\widehat{V}=T_p\widehat{M}$,
$p\in \widehat{M}$)
the corresponding holonomy algebra.
In the present article we shall describe the geometry of the  base $(M,g)$ for each of the three possibilities (i-iii) for the holonomy algebra of the cone $\widehat{M}$.

Our first result describes the holonomy algebra and local structure of a manifold
$(M,g)$ with decomposable holonomy $\hath$ of the cone.

\vskip0.4cm
\noindent
{\bf Theorem \ref{tmain2}.} {\em Let $(M,g)$ be a  pseudo-Riemannian manifold
with decomposable holonomy algebra $\hath$ of the cone $\widehat{M}$.
Then the manifold $(M,g)$ has full holonomy algebra
$\so(p,q)$, where $(p,q)$ is the signature of the metric $g$. Furthermore,
there exists an open dense submanifold $M'\subset M$ such that
any point $p\in M'$ has a neighborhood isometric to a
pseudo-Riemannian manifold of the form
$(a,b)\times N_1 \times N_2$ with the metric given either by
\[g=ds^2+ \cos^2(s) g_1+\sin^2(s)g_2\quad\text{or}\quad g=-ds^2+ \ch^2(s) g_1+\sh^2(s)g_2,\]
where $g_1$ and $g_2$ are metrics on $N_1$ and $N_2$ respectively.}

\vskip0.4cm

Let us recall the following fundamental theorem of Gallot which settles the problem for Riemannian cones over complete Riemannian manifolds.

\bt[S. Gallot, \cite{gallot79}]\label{proposition3.1intro}
Let $(M,g)$ be a complete Riemannian manifold of dimension $\ge2$
with decomposable holonomy
algebra $\hath$ of the cone $\widehat{M}$. Then $(M,g)$
has constant curvature 1 and the cone is flat. 
If, in addition, $(M,g)$ is simply connected, then it is equal to the
standard sphere.
\et

For pseudo-Riemannian manifolds $(M,g)$ the completeness assumption yields
only the following generalisation of Gallot's result:

\vskip0.4cm
\noindent
{\bf Theorem \ref{tmain4}.} {\it \formtmainfore }

\vskip0.4cm

For compact and complete pseudo-Riemannian manifolds $(M,g)$ we are able to establish
the same conclusion as in Theorem \ref{proposition3.1intro}:

\vskip0.4cm
\noindent
{\bf Theorem \ref{tmain3}.} {\it \formtmainthree }

\vskip0.4cm

We remark that for indefinite pseudo-Riemannian manifolds compactness does not imply
completeness, see for example \cite[p. 193]{oneill83} for a geodesically
incomplete Lorentz metric on the 2-torus (the so-called Clifton-Pohl torus).

Since there is no simply connected  compact indefinite
pseudo-Riemannian manifold
of constant curvature 1, we obtain the following
corollary.

\bc
If $(M,g)$ is a simply connected  compact and complete indefinite
pseudo-Riemannian manifold, then the holonomy algebra of the cone
$(\widehat{M},\widehat{g})$ is indecomposable.
\ec

Now we consider the case (ii) when the holonomy algebra $\hath$ of
the cone $\widehat{M}$ is indecomposable but reducible. We completely
analyse the situation in the following two cases:
\begin{itemize}
\item[(ii.a)]  $\widehat{\mathfrak{h}}$
preserves a decomposition
$T_p\widehat{M}=V\+W$ $(p\in \widehat{M})$ into two complementary
subspaces $V$ and $W$.
\item[(ii.b)] $\widehat{M}$ is Lorentzian.
\end{itemize}

In the case (ii.a) one can show that $\widehat{M}$
admits (locally) a para-K\"ahler structure,
which means that the holonomy algebra $\hath$ preserves two complementary
\emph{isotropic} subspaces.  The following theorem characterises
para-K\"ahlerian cones as cones over para-Sasakian manifolds.

\vskip0.4cm
\noindent
{\bf Theorem \ref{tmain6}.} {\it
Let $(M,g)$ be a pseudo-Riemannian manifold. There
is a one-to-one correspondence between para-Sasakian structures
$(M,g,T)$ on $(M,g)$ and para-K\"ahler structures
$(\widehat{M},\widehat{g},\widehat{J})$ on the cone
$(\widehat{M},\widehat{g})$. The correspondence is given by
$T \mapsto \widehat{J} := \widehat{\nabla}T$.
}

\vskip0.4cm

Similarly, we have the following characterisation of the case
when the cone $\widehat{M}$
admits (locally) a para-hyper-K\"ahler structure, which means that the holonomy
algebra $\hath$ preserves two complementary
isotropic subspaces $T^\pm$ and a skew-symmetric
complex structure $J$ such that $JT^+=T^-$. In particular, it preserves
the para-hyper-complex structure $(\widehat{J}_1,\widehat{J}_2,\widehat{J}_3\widehat{J}_1\widehat{J}_2)$, where $\widehat{J}_1|_{T^\pm}=\pm Id$ and
$\widehat{J}_2=J$.

\vskip0.4cm
\noindent
{\bf Theorem \ref{tmain63Sasaki}.} {\it \formtmainsix3Sasaki }

\vskip0.4cm

Finally,  we consider the case (ii.b) when the cone is Lorentzian with
indecomposable
reducible holonomy algebra.

\vskip0.4cm \noindent {\bf Theorem \ref{tmain7}. } {\it Let $(M,g)$ be a Lorentzian manifold of signature $(1,n-1)$ or a
negative definite Riemannian manifold and $(\widehat{M} = \bR^+ \times M, \widehat{g})$ the cone over $M$ with Lorentzian
signature $(1,n)$ or $(n,1)$ respectively. If the holonomy algebra $\hath$ of $\widehat{M}$ is indecomposable reducible
(i.e.\ preserves an isotropic line) then it annihilates a non-zero isotropic vector.}

\vskip0.4cm

The next theorem treats the case of
a  Lorentzian cone $\widehat{M}$ over a negative definite Riemannian
manifold $M$.
\vskip0.4cm

\noindent {\bf Theorem \ref{tmain7}'. } {\it Let $(M,g)$ be a negative definite Riemannian manifold and $(\widehat{M},
\widehat{g})$ the cone over $M$ equipped with the Lorentzian metric of signature  $(+,-,\cdots,-)$. If $\widehat{M}$
admits a non-zero parallel isotropic vector field then $M$ is locally isometric to a manifold of the form
\begin{equation}(M_0=(a,b)\times N,g=-ds^2+e^{-2s}g_N),\label{M0Equ}
\end{equation}
where $a\in\bR\cup\{-\infty\},\,\,b\in\bR\cup\{+\infty\},\,\,a<b$ and
$(N,g_N)$ is a negative definite Riemannian manifold.
Furthermore, if $\hol(\widehat{M_0})$ is indecomposable then
$$\hol(\widehat{M_0},\widehat{g})\cong\hol(N,g_N)\ltimes\bR^{\dim N}.$$ If the manifold $(M,g)$ is complete then the
isometry is global,  $(N,g_N)$ is complete and $(a,b)=\bR$.} \vskip0.4cm

The last theorem treats a  Lorentzian cone $\widehat{M}$ over  Lorentzian manifold $M$.

\vskip0.4cm \noindent {\bf Theorem \ref{tmain7}''. } {\it Let $(M,g)$ be a Lorentzian  manifold and
$(\widehat{M},\widehat{g})$ the cone over $M$ equipped with the Lorentzian metric. If $\widehat{M}$ admits a non-zero
parallel isotropic vector field then there exists an open dense submanifold $M'\subset M$ such that any point of  $M'$ has
a neighborhood isometric to a manifold of the form (\ref{M0Equ}), where $(N,g_N)$ is a positive definite Riemannian
manifold.
Furthermore, if $\hol(\widehat{M_0})$ is indecomposable then
$$\hol(\widehat{M_0},\widehat{g})\cong\hol(N,g_N)\ltimes\bR^{\dim N}.$$
 If  the manifold
$(M,g)$ is complete then each connected component of $M'$ is isometric to a manifold of the form (\ref{M0Equ}), where
$(N,g_N)$ is a positive definite Riemannian manifold and $(a,b)=\bR$.}
\vskip0.4cm
Let us conclude this introduction with some brief remarks about applications of these.
The theorem of Gallot was used by C. B\"{a}r in the classification of Riemannian manifolds admitting a real Killing spinor \cite{baer93}. In general, a pseudo-Riemannian manifold  admits a real/imaginary Killing spinors if and only if its space-like/time-like cone admits a parallel spinor (details in \cite{Bohle}). Hence, a strategy for studying  manifolds with Killing spinors is to study their cones with parallel spinors. Now, in order to classify manifolds with parallel spinors the knowledge of their holonomy group is essential.  In the Riemannian situation Gallot's result reduces the problem to irreducible holonomy groups of cones. With our results this strategy becomes applicable to arbitrary signature.

Another applications in the same spirit --- solutions to an overdetermined system of PDE's correspond to parallel sections for a certain connection --- comes from conformal geometry. Here, to a conformal class of a metric one can assign the so-called Tractor bundle with Tractor connection. Parallel sections of this connection correspond to metrics in the conformal class which are Einstein. For conformal classes which contain an Einstein metric the holonomy of the conformal Tractor connection reduces to the Levi-Civita holonomy of the Fefferman-Graham ambient metric \cite{fefferman/graham85}.  For conformal classes containing proper Einstein metrics with positive/negative Einstein constant the ambient metric reduces to the space-like/time-like cone over a metric in the conformal class \cite{leitner05, armstrong05,armstrong-leistner06}. 
Again, our results enable us to describe the holonomy of the conformal Tractor connection by the holonomy of the cone.

To carry out the details of both applications lies beyond the scope of this paper and will be subject to future research.

\paragraph{Acknowledgements.} The authors thank the International Erwin Schr\"odinger Institute for Mathematical Physics in Vienna for the hospitality during the Special Research Semester {\em Geometry of Pseudo-Riemannian Manifolds with Applications in Physics}.
\section{Doubly warped products}
Let $(N_1,g_1)$ and $(N_2,g_2)$ be two pseudo-Riemannian manifolds, $\e =\pm 1$
and $f_1$, $f_2$ two nowhere vanishing
smooth functions on an open interval $I=(a,b)\subset \bR$.
The manifold $M=I\times N_1\times N_2$ with the metric
\begin{equation}g=\e ds^2+f_1(s)^2g_1 +f_2(s)^2g_2 \label{dwEqu}\end{equation}
is called a \emph{doubly warped product}. We will consider the
coordinate vector field $\partial_0=\partial_s$ on the interval $I$ and
vector fields $X, X',\ldots $ on $N_1$ and $Y, Y',\ldots $ on $N_2$ as
vector fields on $M$.
\bs \label{doubleProp}
\begin{itemize}
\item[(i)] The Levi-Civita connection $\n$ of the pseudo-Riemannian manifold
$(M,g)$ is given by:
\begin{eqnarray*}
\n_{\partial_0}\partial_0&=&0\\
\n_{\partial_0}X&=&\n_X\partial_0=\frac{f_1'}{f_1}X\\
\n_{\partial_0}Y&=&\n_Y\partial_0=\frac{f_2'}{f_2}Y\\
\n_XX'&=&-\e\frac{f_1'}{f_1}g(X,X')\partial_0+\n^1_XX'\\
\n_YY'&=&-\e\frac{f_2'}{f_2}g(Y,Y')\partial_0+\n^2_YY'\\
\n_XY&=&\n_YX=0,
\end{eqnarray*}
where $\n^1$ and $\n^2$ are the Levi-Civita connections of
$g_1$ and $g_2$.
\item[(ii)] A curve $I\ni t\mapsto (s(t),\g_1(t),\g_2(t))\in
I\times N_1\times N_2=M$
is a geodesic if and only if it satisfies the equations:
\begin{eqnarray*}
\ddot{s} &=& \e (f_1'(s)f_1(s)g_1(\dot{\g_1},\dot{\g_1})+
f_2'(s)f_2(s)g_2(\dot{\g_2},\dot{\g_2}))\\
\n^i_{\dot{\g_i}}\dot{\g_i} &=& -2\frac{f'_i}{f_i}\dot{s}\dot{\g_i},\; i=1,2.
\end{eqnarray*}
\item[(iii)] In terms of the arclength parameters $u_1$, $u_2$ of
the curves $\g_1$, $\g_2$ in $N_1$, $N_2$ the equations (ii) take the
form:
\begin{eqnarray} \label{geoEqu}
\ddot{s} &=& \e\e_1 f_1'(s)f_1(s)\dot{u_1}^2+
\e \e_2f_2'(s)f_2(s)\dot{u_2}^2\\ \label{geoEqu2}
\ddot{u_i} &=& -2\frac{f'_i}{f_i}\dot{s}\dot{u_i},\; i=1,2.
\end{eqnarray}
where $\e_i = g_i(\frac{d\g_i}{du},\frac{d\g_i}{du}) \in \{ \pm 1, 0\}$.
\end{itemize}
\es

\bfolg \begin{itemize}
\item[(i)] The submanifolds $I\times N_i \subset M$, $i=1, 2$,
are totally geodesic, as well as their intersection $I$.
\item[(ii)] Any geodesic $\G\subset M$ is contained in a totally geodesic
submanifold  $I\times \G_1\times \G_2 \subset M$, where
$\G_i={pr}_{N_i}(\G )\subset N_i$ are geodesics in $N_i$, $i=1,2$.
\end{itemize}
\efolg

\bs The curvature tensor of the doubly warped product (\ref{dwEqu})
is given by:
\begin{eqnarray*}
R(\partial_0,X)&=&-\e \frac{f_1''}{f_1} \partial_0\wedge X\\
R(\partial_0,Y)&=&-\e \frac{f_2''}{f_2} \partial_0\wedge Y\\
R(X,X')&=&-\e \left( \frac{f_1'}{f_1}\right)^2 X\wedge X' + R^1(X,X')\\
R(Y,Y')&=&-\e \left( \frac{f_2'}{f_2}\right)^2 Y\wedge Y' + R^2(Y,Y')\\
R(X,Y)&=& -\e \frac{f_1'f_2'}{f_1f_2}X\wedge Y
\end{eqnarray*}
\es Recall that a pseudo-Riemannian metric $g$ has constant curvature $\kappa$  if and only if its curvature tensor has
the form $$R(X,Y) = \kappa X\wedge Y =\kappa (X\ot g(Y,\cdot )- Y\ot g(X,\cdot )).$$ \bfolg A doubly warped product
(\ref{dwEqu}) has constant curvature $\kappa=-\e c$ if and only if the metrics $g_1$ and $g_2$ have constant curvature
$\kappa_1$, $\kappa_2$ and the warping functions satisfy the following system of equations: $$\frac{f_1''}{f_1}= c \quad
\mbox{if}\quad \dim N_1>0,$$ $$\frac{f_2''}{f_2}= c\quad \mbox{if}\quad \dim N_2>0,$$ $$\frac{f_1'f_2'}{f_1f_2}=c,  \quad
\mbox{if}\quad \dim N_1>0\quad\mbox{and}\quad \dim N_2>0,$$ $$\left( \frac{f_1'}{f_1}\right)^2 -\e
\frac{\kappa_1}{f_1^2}=c\quad \mbox{if}\quad \dim N_1>1,$$ $$\left( \frac{f_2'}{f_2}\right)^2 -\e
\frac{\kappa_2}{f_2^2}=c\quad \mbox{if}\quad \dim N_2>1.$$ \efolg

\noindent
Solving these equations we get the following corollary. We will denote
by $g_k, g'_k,g''_k$ pseudo-Riemannian metrocs of constant
curvature $k \in \{ \pm 1 , 0\}$.
\bfolg \label{ccCor}  Then the following
doubly warped product metrics $g_k$ have constant curvature $k$:
\begin{eqnarray*}
g_{-\e}&=&\e ds^2 + \cosh^2(s)g_{-\e }'+\sinh^2(s)g_\e''\\
g_{\e}&=&\e ds^2 + \cos^2(s)g_{\e }'+\sin^2(s)g_\e''\\
g_{-\e}&=&\e ds^2 + e^{2s}g_0'\\
g_0&=&\e ds^2 +s^2g_\e'+g_0''\\
g_{-\e}&=&\e ds^2 + \cosh^2(s)dt^2+\sinh^2(s)g_\e''\\
g_{-\e}&=&\e ds^2 + \cosh^2(s)g_{-\e }'+\sinh^2(s)du^2\\
g_{\e}&=&\e ds^2 + \cos^2(s)dt^2+\sin^2(s)g_\e''\\
g_{\e}&=&\e ds^2 + \cos^2(s)g_{\e }'+\sin^2(s)du^2\\
g_0&=&\e ds^2 +s^2dt^2+g_0''\\
g_{-\e}&=&\e ds^2 +\sinh^2(s)g_\e'\\
g_{-\e}&=&\e ds^2 + \cosh^2(s)g_{-\e }'\\
g_{\e}&=&\e ds^2 + \sin^2(s)g_\e'\\
g_0&=&\e ds^2 +g_0'\\
g_{-\e}&=&\e ds^2 \pm \cosh^2(s)dt^2\pm \sinh^2(s)du^2\\
g_{\e}&=&\e ds^2 \pm \cos^2(s)dt^2\pm \sin^2(s)du^2\\
g_{-\e}&=&\e ds^2 \pm \sinh^2(s)dt^2\\
g_{-\e}&=&\e ds^2 \pm \cosh^2(s)dt^2\\
g_{\e}&=&\e ds^2 \pm \sin^2(s)dt^2\\
\end{eqnarray*}
Any doubly warped product of pseudo-Riemannian manifolds
$(N_1,g_1)$, $(N_2,g_2)$ which has constant curvature $\pm 1$ or $0$
belongs to the above list up to a shift $s\mapsto s+s_0$ and
rescaling $$(f_i^2,g_i)\mapsto (\lambda_i^2f_i^2,\frac{1}{\lambda_i^2}g_i).$$
\efolg
\subsection*{Geometric realisation of doubly warped products
of constant curvature}
Now we give a realisation of the above doubly warped
products in terms the pseudo-sphere models of the spaces
of constant curvature.
\subsubsection*{ The standard pseudo-spheres as models of
spaces of curvature $\pm 1$}
Let $\bR^{t,s}=(\bR^{t+s},\langle \cdot ,\cdot \rangle -\sum_{i=1}^tdx_i^2+\sum_{i=t+1}^{t+s}dx_i^2)$ be the
standard pseudo-Euclidian vector space of signature $(t,s)$.
We denote by
\begin{eqnarray*}
S^{t,s}_+&:=&\{ x\in \bR^{t,s+1}| \langle x, x\rangle  = +1\}\\
S^{t,s}_-&:=&\{ x\in \bR^{t+1,s}| \langle x, x\rangle  = -1\}\\
\end{eqnarray*}
the two unit pseudo-spheres.  The induced metric $g_\pm = g_{S^{t,s}_\pm}$
of $S^{t,s}_\pm$ has
signature $(t,s)$ and constant curvature $\pm 1$. More precisely
the curvature tensor is given by
\[ R(X,Y)Z = \pm (\langle Y,Z\rangle X -\langle X,Z\rangle Y).\]
Notice that $S^{0,n}_+=S^n$ is the standard unit n-sphere,
$S^{0,n}_-=H^n$ is hyperbolic n-space, $S^{1,n-1}_+=dS^n$ is de Sitter
n-space and  $S^{1,n-1}=AdS^n$ is anti de  Sitter n-space.

\subsubsection*{Flat space as cone over the pseudo-spheres}
The domains
$$\bR^{t,s}_\pm :=\{ x\in \bR^{t,s}| \pm \langle x, x\rangle  >0\}\subset
\bR^{t,s}$$
are isometrically identified via the map $(r,x)\mapsto rx$
with the  space-like or time-like cone over $S^{t,s-1}_+$
or $S^{t-1,s}_-$ endowed with the metric $\pm dr^2 + r^2g_\pm$, respectively.
In particular, the space-like cone over a space of constant curvature $1$
and the time-like cone over a space of constant curvature $-1$
are flat.

\subsubsection*{Realisation of doubly warped products by double
polar coordinates}
Now we show that any splitting of a pseudo-Euclidian vector space
as an orthogonal sum of two pseudo-Euclidian subspaces
induces local  parametrisations of the
pseudo-spheres. Using these 'double
polar' parametrisations (more precisely, polar equidistant
parametrisations \cite{alekseevskij-vinberg-solodovnikov93})
we will  show that the spaces of constant
curvature can be locally presented as
doubly warped products with trigonometric or hyperbolic warping functions
over spaces of appropriate  constant curvature.

We consider the pseudo-spheres $S_+(V)=S_+^{t,s}\subset V=\bR^{t,s+1}$ and
$S_-(V)=S_-^{t,s}\subset V=\bR^{t+1,s}$. Any orthogonal decomposition
$V=V_1\+V_2 =\{ v=x+y| x\in V_1, y\in V_2\}$, $
\langle\cdot ,\cdot \rangle=\langle\cdot ,\cdot \rangle_1
+\langle\cdot ,\cdot \rangle_2$ defines  a diffeomorphism
\begin{equation} \label{DEqu}
(s,\bar{x},\bar{y}) \mapsto x+y,\quad x=\cos (s)\bar{x},
\quad y=\sin (s)\bar{y},
\end{equation}
of
$(0,\frac{\pi}{2})\times S_\e (V_1)\times S_\e (V_2)$ onto the (not necessarily
connected) domain
$$D = \{ v=x+y\in S_\e(V)| 0<\e \langle x, x\rangle_1
<1\} .$$
Similarly the map
\begin{equation} \label{D'Equ}(s,\bar{x},\bar{y}) \mapsto x+y,\quad
x=\cosh (s)\bar{x},\quad y\sinh (s)\bar{y},
\end{equation}
is a diffeomorphism of $\bR^+  \times S_\e (V_1)\times S_{-\e} (V_2)$
onto the  domain
$$D' = \{ v=x+y\in S_\e(V)| \e \langle x, x\rangle_1 >1
\} .$$
\bs With respect to the diffeomorphisms
(\ref{DEqu}) and (\ref{D'Equ})  the metric $g_\e$ of $S_\e (V)$
is given by
\begin{eqnarray*}
g_\e|_{D} &=&  \e ds^2 + \cos^2(s)g_{S_\e(V_1)}+\sin^2(s)g_{S_\e(V_2)}\\
g_\e|_{D'} &=&  -\e ds^2 + \cosh^2(s)g_{S_\e(V_1)}+\sinh^2(s)g_{S_{-\e}(V_2)}.
\end{eqnarray*}
\es
\subsubsection*{Horospherical coordinates and corresponding
warped products}
Let $(V,\langle \cdot ,\cdot \rangle )$ be an indefinite pseudo-Euclidian
vector space, $p, q\in V$ two isotropic vectors such that $\langle p,q\rangle
=1$ and $W=\mbox{{\rm span}}\{ p,q\}^\perp$. Then
\[ \bR^+\times W \ni (s,\xi )\mapsto y=up+vq+x\in S_\e (V),\quad u=\pm
\frac{1}{2}e^{-s}(\e -e^{2s}
\langle \xi ,\xi \rangle ),\quad v=\pm 2e^s,\quad x=e^s\xi\]
is a diffeomorphism onto the domain $S_\e (V)\cap \{ y\in V|\pm v>0\}$.
In the coordinates $(s,\xi )$ the hypersurfaces $s=const$ correspond
to the  hyperplane sections (horospheres) $\{y\in S_\e(V)
| \langle y,p\rangle = \pm e^s\}$ and the curves $\xi =\xi_0=const$
are geodesics perpendicular to the horospheres.
A direct calculation shows that:
\bs  \label{horoProp} The induced metric of the pseudo-sphere $S_\e (V)$ in horospherical
coordinates $(s,\xi )$ is given by:
$$g_\e=\e ds^2+e^{2s}g_0,$$
where $g_0=d\xi^2$ is the induced pseudo-Euclidian metric on $W$.
\es
\subsubsection*{Completeness of some doubly warped products}
\bs \label{completeProp}
Let $(N_1,g_1)$, $(N_2,g_2)$ be pseudo-Riemannian manifolds
and $(M=I\times N_1\times N_2, g=\e ds^2 + f_1(s)^2g_1+f_2(s)^2g_2)$
a doubly warped product with non-constant
warping functions as in Corollary \ref{ccCor}.
Then $(M,g)$ is complete only in the following cases:
\begin{itemize}
\item[(i)]
\[ g=\e ds^2 + \cosh^2 (s)g_1,\]
where $I=\bR$ and $g_1$ is complete,  and
\item[(ii)]
\[ g=\e ds^2 + e^{2s}g_1\]
where $I=\bR$ and $\e g_1$ is  complete and positive definite.
\end{itemize}
\es

\bprf The system
(\ref{geoEqu}-\ref{geoEqu2}) has solutions
$u_i=u_i^0=const$, $s=at+b$, which are
complete if and only if $I=\bR$. This excludes all the warping functions
which have a zero. It remains to check that the metric (i) is
complete for any complete metric $g_1$  and that (ii) is complete only
if $\e g_1$ is complete and positive definite.
In fact, in both cases the squared velocity $l= g(\dot{\g},\dot{\g})$ is
constant. In the second case, for instance, it is given by
$l=\e \dot{s}^2 +\e_1
e^{2s}\dot{u}^2$, $u:=u_1$, which yields $\ddot{y}=\e ly$ after
the substitution $y=e^s$. The differential equation  $\ddot{y}=\e ly$
admits solutions which are positive on the
real line if and only if $\e l>0$. The positivity is necessary since
$y=e^s$. This shows that $g$ is positive
or negative definite, i.e.\ $\e g_1$ is positive definite.  The other case
is similar, see \cite{Bohle}, where the case of Lorentzian signature
is considered.
\eprf

\section{Examples of cones with reducible holonomy}
Let $\widehat{g}=c dr^2 + r^2g$ be the cone metric on
$\widehat{M}:=\rr^+\times M$, where $(M,g)$ is a pseudo-Riemannian manifold.
Depending on the sign of the constant $c$ the cone is called space-like
$(c>0)$ or time-like $(c<0)$.
Later on we will assume, without restriction of generality, that $c=1$.
In fact,  as we allow $g$ to be
of any signature we can rescale $\widehat{g}$ by
$\frac{1}{c}\in \bR^*$.

We denote by $\p_r$ the radial vector field.
The  Levi-Civita connection of the cone $(\widehat{M},\widehat{g})$ is given by
\beq
\left.\begin{array}{rcl}
\widehat{\n}_{\p_r}\p_r  &=& 0,\\
\widehat{\n}_X\p_r  &=& \frac{1}{r} X,\\
\widehat{\n}_X Y &=& \n_X Y -\frac{r}{c} g(X,Y)\p_r,
\end{array}\right\}\label{lem1}
\eeq
for all vector fields $X,Y\in \Gamma(T\widehat{M})$ orthogonal to $\p_r$.
The curvature $\hut{R}$ of the cone is given by the following formulas
including the curvature $R$ of the base metric $g$:
\beq
\left.\begin{array}{rcl}
\p_r\inter \hut{R}&=&0,\\
\hut{R}(X,Y)Z&=& R (X,Y)Z - \frac{1}{c}\left(  g(Y,Z)X - g(X,Z)Y \right), \text{ or}\\
\hut{R}(X,Y,Z,U)&=& r^2\left( R (X,Y,Z,U) - \frac{1}{c}\left(g(Y,Z)g(X,U) - g(X,Z)g(Y,U)  \right)\right),
\end{array}\right\}\label{lem2}
\eeq
for $X,Y,Z,\ U\in TM$.
This implies that if $(M,g)$ is a space of constant curvature $\kappa$, i.e.
\be
 R(X,Y,Z,U) &=& \kappa\left(g(X,U)g(Y,Z) - g(X,Z)g(Y,U)\right),
\ee
the cone has the curvature $r^2\left(\kappa-\frac{1}{c}\right) \left(g(X,U)g(Y,Z) - g(X,Z)g(Y,U)\right)$.
In particular, if $\kappa=\frac{1}{c}$, then the cone is flat,
as it is the case for the $c=1$ cone over the standard sphere of radius
$1$ or the $c=-1$ cone over the hyperbolic space.

{}From now on we assume $c=\pm 1$.
 We denote by $(\widehat{M}=\bR^+\times M, \widehat{g}=dr^2+r^2g)$, the
space-like cone over $(M,g)$ and by $(\widehat{M}^-=\bR^+\times M,
\widehat{g}^-=-dr^2+r^2g)$ the time-like cone.
Notice that
the metric $\widehat{g}^-$ of the time-like cone
$\widehat{M}^-$ over $(M,g)$ is obtained
by multiplying the metric $dr^2-r^2g$ of the
space-like cone  over $(M,-g)$ by $-1$. Thus it is sufficient
to consider only  space-like cones.

We will now present some examples which illustrate that Gallot's statement is not
true in arbitrary signature, and that the assumption of completeness is essential even in the Riemannian situation.

\bbsp\label{ex1} Let $(F,g_F)$ be a complete pseudo-Riemannian manifold of
dimension at least $2$ and which is not of constant curvature $1$.
Then the pseudo-Riemannian manifold $$(M=\bR\times F,g=-ds^2+\ch^2(s)g_F)$$
is complete, the restricted  holonomy group of the cone over $(M,g)$ is
non-trivial and admits a non-degenerate invariant proper subspace.
\ebsp
\bprf
The manifold $(M,g)$ is complete by Proposition
\ref{completeProp}.
The  non-vanishing terms of the Levi-Civita connection $\n$ of $(M,g)$
are given by
\beq
\left.\begin{array}{rcl}
{\n}_X \p_s &=& \th(s) X,\\
{\n}_{\p_s} X &=& \p_s X+\th(s) X,\\
{\n}_X Y &=& \n^F_X Y +\ch(s)\sh(s) g_F(X,Y)\p_s,
\end{array}\right\}
\eeq where  $X,Y\in TF$ are vector fields depending on the parameter $s$ and $\n^F$ is the Levi-Civita connection of
the manifold $(F,g_F)$.
Consider on $\widehat{M}$ the vector field $X_1=\ch^2(s)\p_r-\frac{1}{r}\sh(s)\ch(s)\p_s$.
We have $\widehat{g}(X_1,X_1)=\ch^2(s)>0$. It is easy to check that the distribution generated by
the vector field $X_1$ and by the distribution $TF\subset T\widehat{M}$ is parallel.

For the  curvature tensor $R$ of $(M,g)$ we have
\beq
\begin{array}{rcl}
 R(X,Y)Z&=&  R_F (X,Y)Z + \tanh^2(s)\left(g_F(Y,Z)X - g_F(X,Z)Y \right),
\end{array}
\eeq
where $X,Y,Z,\ U\in TF$ and $R_F$ is the curvature tensor  of $(F,g_F)$.
This shows that $(M,g)$ cannot have constant sectional curvature, unless
$F$ has constant curvature $1$ (see Corollary \ref{ccCor}).
Thus the cone
$(\widehat{M},\widehat{g})$ is not flat. \eprf

\bbsp  \label{horosphere} Let $M$ be a manifold
of the form $\bR\times N$
with the metric $g=-(dt^2+e^{-2t}g_N)$, where $(N,g_N)$ is a
pseudo-Riemannian manifold.
Then \begin{itemize} \item[1.] The light-like vector field
$e^{-t}(\p_r+\frac{1}{r}\p_t)$ on the space-like cone $\widehat{M}$ is parallel.
\item[2.] If $(N=N_1\times N_2,g_N=g_1+g_2)$ is a product of a flat
pseudo-Riemannian manifold $(N_1,g_1)$ and of a non-flat
pseudo-Riemannian manifold $(N_2,g_2)$, then $\widehat{M}$ is
locally decomposable and not flat\footnote{We learned this from Helga Baum}.
In fact, there is a parallel
non-degenerate  \emph{flat} distribution of dimension $\dim N_1$
on $\widehat{M}$. \end{itemize}
\ebsp

The manifold $(M,g)$ in Example \ref{horosphere} is complete if and only $g_N$ is complete and positive definite, see
Proposition \ref{completeProp}. Notice that $g$ is the hyperbolic metric in horospherical coordinates if $(N,g_N)$ is
Euclidian space.

\bbsp \label{ex2} Let $(M_1,g_1)$ and $(M_2,g_2)$ be two pseudo-Riemannian manifolds.
Then the product of the cones $$(\widehat{M_1}\times \widehat{M_2}= (\bR^+ \times M_1)\times(\bR^+ \times M_2), \widehat{g}=(dr_1^2 + r_1^2 g_1)+(dr_2^2 + r_2^2 g_2))$$
is the space-like cone over the manifold
$$(M=\left(0,\frac{\pi}{2}\right)\times M_1\times M_2,g=ds^2+\cos^2(s)g_1+\sin^2(s)g_2).$$\ebsp
\bprf Consider the functions $$r=\sqrt{r_1^2+r_2^2}\in\bR^+,\quad s=\arctg\left(\frac{r_2}{r_1}\right)\in \left(0,\frac{\pi}{2}\right).$$
Since $r_1,r_2>0$, the  functions  $r$ and $s$ give a diffeomorphism
$\bR^+\times\bR^+\to \bR^+\times \left(0,\frac{\pi}{2}\right)$.
For   $\widehat{M_1}\times \widehat{M_2}$ we get
$$\widehat{M_1}\times \widehat{M_2}\cong \bR^+\times \left(0,\frac{\pi}{2}\right)\times M_1\times M_2$$
and
$$\widehat{g_1}+\widehat{g_2}=dr^2+r^2 (ds^2+\cos^2(s)g_1+\sin^2(s)g_2).$$ \eprf

Suppose that the manifolds $(M_1,g_1)$ and $(M_2,g_2)$ are Riemannian.
Then  the manifold $(M,g)$ is Riemannian and incomplete. The cone
over $M$ is decomposable. Moreover, it is not flat, unless the manifolds
$(M_1,g_1)$, $(M_2,g_2)$ are of dimension less than $2$
or of constant curvature $1$, see Corollary \ref{ccCor}.
Example \ref{ex2} shows that the completeness assumption in Theorems
\ref{proposition3.1intro} and \ref{tmain3}
is necessary.

\bbsp \label{ex3} Let $(M_1,g_1)$ and $(M_2,g_2)$ be two pseudo-Riemannian
manifolds. Then the space-like cone over the manifold
$$(M=\bR^+\times M_1\times M_2,g=-ds^2+\ch^2(s)g_1+\sh^2(s)g_2)$$
is isometric to the open subset $\Omega =\{ r_1 >r_2\}$ in the product of the
cones
$$(\widehat{M_1}\times \widehat{M}_2
= (\bR^+ \times M_1)\times(\bR^+ \times M_2), \widehat{g_1}+\widehat{g_2}^-=(dr_1^2 + r_1^2 g_1)+(-dr_2^2 + r_2^2 g_2)).$$
\ebsp
\bprf Consider the functions
$$ r=\sqrt{r_1^2-r_2^2}\in\bR^+,\quad s= \arcth \left(\frac{r_2}{r_1}\right)\in\bR^+ .$$
The functions  $r$ and $s$ give a diffeomorphism $\{ (r_1,r_2)\in \bR^2|
0<r_2<r_1\} \to \bR^+\times \bR^+$.
For   $\widehat{M_1}\times \widehat{M_2}$ we get
$$\Omega \cong \bR^+\times \bR^+\times M_1\times M_2$$
and
$$\widehat{g_1}+\widehat{g_2}^-=dr^2+r^2 (-ds^2+\ch^2(s)g_1+\sh^2(s)g_2).$$ \eprf

\bbsp \label{ex4}
Let $(t,x_1, \ldots, x_n, x_{n+1}, \ldots , x_{2n}) $ be coordinates on $\rr^{2n+1}$. Consider  the
metric $g$ given by
\be
 g&=& \left(
 \begin{array}{c|c|c}
 -1&0&u^t \\
 \hline
 0&0&H^t \\
 \hline
 u & H & G
 \end{array}\right)
 \ee
 where
 \bi
 \item $u=(u_1,\ldots, u_n)$ is a diffeomorphism of $\rrn$, depending on $x_1,\ldots , x_n$,
 \item  $H = \einhalb\left(\ddxj(u_i)\right)_{i,j=1}^n$ its  non-degenerate Jacobian,
  and
  \item $G$ the symmetric matrix given by $G_{ij}= -u_iu_j$.
  \ei Then  the space-like cone over $(\rr^{2n+1},g)$ is not flat but its holonomy representation decomposes into two totally isotropic invariant subspaces.
   For the proof of this see Proposition \ref{psbsp} in  Section \ref{isosec}.
 \ebsp

\section{Local structure of decomposable cones}
\label{LocalSec}

In this section we assume that the holonomy group of the cone $(\widehat{M},\widehat{g})$ is decomposable and we give a local description of the manifold $(M,g)$, independently of completeness.

Suppose that
the  holonomy group $Hol_x$ of $(\widehat{M},\widehat{g})$ at a point $x\in\widehat{M}$ is decomposable, that is $T_x\widehat{M}$
is a sum $T_x\widehat{M}=(V_1)_x\oplus (V_2)_x$ of two non-degenerate $Hol_x$-invariant orthogonal subspaces.
They define two parallel non-degenerate distributions $V_1$ and $V_2$.
Denote by  $X_1$ and $X_2$ the projections of the vector field $\p_r$ to the distributions $V_1$ and $V_2$ respectively. We have
\beq\label{e1} \p_r=X_1+X_2.\eeq
We decompose the vectors $X_1$ and $X_2$ with respect to the decomposition $T\widehat{M}=T\bR^+ \oplus TM$,
\beq\label{e2} X_1=\a\p_r+X,\quad  X_2=(1-\a)\p_r-X,\eeq
where $\a$ is a function on $\widehat{M}$ and $X$ is a vector field on $\widehat{M}$ tangent to $M$.
We have \beq\label{e3}\widehat{g}(X,X)=\a-\a^2,\quad \widehat{g}(X_1,X_1)=\a, \quad \widehat{g}(X_2,X_2)=1-\a.\eeq

\blem \label{lemA3} The open subset $U=\{x|\a(x)\neq 0,1 \}\subset\widehat{M}$ is dense.\elem
\bprf  Suppose that $\a=1$ on an open subspace  $V\subset\widehat{M}$. We claim that  $\p_r\in V_1$ on $V$.
Indeed, on $V$  we have $$X_1=\p_r+X,\,\,\, X_2=-X\,\,\,\,\text{ and   }\,\,\,\,\, \widehat{g}(X,X)=0.$$
We show that $X=0$.
Let $Y_2\in V_2$. We have the decomposition $Y_2=\l\p_r+Y,$ where $\l$ is a function on $\widehat{M}$ and $Y\in TM$.
It is  $$\widehat{\n}_{Y_2}X_1=\widehat{\n}_{\l\p_r+Y}(\p_r+X)=\frac{1}{r}Y+\widehat{\n}_{Y_2}X.$$
Note that $Y=Y_2-\l X_1+\l X$. Hence, $$\widehat{\n}_{Y_2}X_1=\frac{1}{r}(Y_2-\l X_1+\l X)+\widehat{\n}_{Y_2}X.$$
Since $X,\widehat{\n}_{Y_2}X,Y_2\in V_2$ and $\widehat{\n}_{Y_2}X_1\in V_1$, we see that $$\frac{1}{r}(Y_2+\l X)+\widehat{\n}_{Y_2}X=0.$$
From $\widehat{g}(X,X)=0$ it follows that $\widehat{g}(\widehat{\n}_{Y_2}X,X)=0$. Thus we get $\widehat{g} (Y_2,X)=0$ for all $Y_2\in V_2$.
Since $V_2$ is non-degenerate, we conclude that $X=0$. Thus
$\partial_r\in V_1$.

Let $Y_2\in  V_2$, then $\widehat{\n}_{Y_2}\p_r=\frac{1}{r}Y_2$.
Since the distribution $V_1$ is parallel and $\p_r\in V_1$, we see that $Y_2=0$ and $V_2=0$.
Contradiction.
\eprf

We now consider the dense open submanifold $U\subset\widehat{M}$.
The vector fields $X_1$, $X_2$ and $X$  are nowhere isotropic on $U$.
For $i=1,2$ let $E_i\subset V_i$ be the subdistribution of $V_i$ orthogonal to $X_i$. Denote by $L$ the distribution of lines on $U$ generated
by the vector field $X$.
We get on $U$ the orthogonal decomposition $$T\widehat{M}=T\bR\oplus L\oplus E_1\oplus E_2.$$

\blem \label{lemA4} Let $Y_1\in  E_1$ and $Y_2\in  E_2$, then  on $U$ we have
\begin{itemize}
\item[1.] $Y_1\a=Y_2\a=\p_r\a=0.$
\item[2.] $\widehat{\n}_{Y_1}X=\frac{1-\a}{r}Y_1$, $\widehat{\n}_{Y_2}X=-\frac{\a}{r}Y_2$.
\item[3.] $\widehat{\n}_{\p r}X=\p_r X+\frac{1}{r}X=0$.
\item[4.] $\widehat{\n}_{X}X=\left(\frac{(1-\a)^2}{r}-X\a \right)X_1+
\left(\frac{\a^2}{r}-X\a \right)X_2$.
\end{itemize}\elem

\bprf Using (\ref{e2}), we have $$\begin{array}{l}
\widehat{\n}_{Y_1}X_1=(Y_1\a)\p_r+\frac{\a}{r}Y_1+\n_{Y_1}X,\\
\widehat{\n}_{Y_1}X_2=-(Y_1\a)\p_r+\frac{1-\a}{r}Y_1-\n_{Y_1}X.\end{array}$$
Since $Y_1\in  E_1\subset  V_1$ and  the distributions $V_1$, $V_2$ are
parallel, projecting these equations  onto $V_2$ and adding them
yields $\widehat{\n}_{Y_1}X_2=0$. Then, from the second equation,
we see that $Y_1\a=0$ and $\widehat{\n}_{Y_1}X=\n_{Y_1}X=\frac{1-\a}{r}Y_1$.
The other claims can be proved similarly. \eprf

\noindent
Since $\p_r\a=0$, the function $\a$ is a function on $M$.
Note that  $U=\bR^+\times U_1$,
where $$U_1=\{x\in M|\a(x)\neq 0,1 \}\subset M.$$

\noindent
Claim 3 of Lemma \ref{lemA4} shows that $X= \frac{1}{r}\t X$, where $\t X$ is a vector field on the manifold $M$.
Hence the distributions $L$ and $E=E_1\oplus E_2$ do not depend on $r$ and  can be considered
as distributions on $M$. Claim 2 of Lemma \ref{lemA4} shows that the distributions $E_1$ and $E_2$ also
do not depend on $r$. We get on $U_1$  the orthogonal decompositions $$TM=L\oplus E,\quad E=E_1\oplus E_2.$$

\blem \label{lemA5} The function $\a$ satisfies on $U_1$ the following
differential equation
     $$\t X\a=2(\a-\a^2).$$\elem
\bprf {}From
\[ r^2\widehat{\n }_XX= \widehat{\n }_{\tilde{X}}\tilde{X} =
\n_{\tilde{X}}\tilde{X}-rg(\tilde{X},\tilde{X})\partial_r =
\n_{\tilde{X}}\tilde{X}-r(\a -\a^2)\partial_r\]
and  Claim 4 of Lemma \ref{lemA4} we conclude that
$\n_{\tilde{X}}\tilde{X}\in TM$ is a linear combination
of $X_1$ and $X_2$ and hence proportional to $X=(1-\a )X_1-\a X_2$.
This implies
\[ 
(2\a-1)\left(X\a-\frac{2}{r}(\a-\a^2)\right)=0.\]
If $\a=\frac{1}{2}$, then $\widehat{\n}_{X}X=\frac{1}{4r}\p_r$ and $\n_{\t X}\t X=\frac{r}{2}\p_r$. The last equality is impossible. \eprf

\bfolg \label{gradCor} On $M$ we have
$$\tilde{X}=\frac{1}{2}{grad}(\a ).$$
\efolg

\noindent
Lemma  \ref{lemA5} implies that if $t$ is a coordinate on $M$ corresponding
to the vector field $\t X$, then
$$\a(t)=\frac{e^{2t}}{e^{2t}+c},$$ where $c$ is a constant.

\noindent
From Lemmas \ref{lemA4} and \ref{lemA5} it follows that
$$\widehat{\n}_{X}X=-\frac{\a-\a^2}{r}\p_r+\frac{1-2\a}{r}X.$$
On the subset $U_1\subset M$ we get the following
\beq\label{e4}  \n_{\t X}\t X=(1-2\a)\t X,\quad \n_{Y_1}\t X=(1-\a)Y_1,\quad \n_{Y_2}\t X=-\a Y_2,\eeq
for  any  $Y_1\in \Gamma (E_1)$ and $Y_2\in \Gamma (E_2)$.

\btheo \label{tmain2} \formtmaintwo \etheo

\bprf Since the distributions $V_1$ and $V_2$ are parallel, for any $Y_1\in V_1$ and $Y_2\in  V_2$ we have
$\widehat{R}(Y_1,Y_2)=0$. From  (\ref{lem2}) it follows that $\widehat{R}(X,Y_2)=\widehat{R}(X,Y_1)=0$. Hence,
$R(\t X,Y)=\t X\wedge Y$
for all vector fields $Y$  on $M$. By Lemma \ref{lemA3}, there exists $y\in M$ such that $g_y(\t X,\t X)\neq 0$.
The holonomy algebra of the manifold $M$ at the point $y$ contains the subspace $\t X_y\wedge T_yM$.
Since $g_y(\t X,\t X)\neq 0$, this vector subspace generates the whole Lie algebra $\so(T_yM,g_y)$. \eprf

\btheo \label{tmain5}
\formtmainfive
\etheo

\bprf We need the following
\blem \label{lemA6} \begin{enumerate}
\item[(i)]
The distributions $E_1, E_2, E=E_1\oplus E_2
\subset TM$ defined on
$U_1\subset M$ are involutive and the distributions $E_1\oplus L, E_2\oplus L\subset TM$ are parallel on $U_1$.
\item[(ii)]
Let $x\in U_1$ and $M_x\subset U_1$ the
maximal connected integral submanifold of the distribution $E$. Then
the distributions $E_1|_{M_x}, E_2|_{M_x}\subset TM_x=E|_{M_x}$ are
parallel.
\end{enumerate}
\elem  \bprf
(i) On $U\subset \widehat{M}$, the distribution $E_i=V_i\cap TM$ is
the intersection of two
involutive distributions and hence involutive, for $i=1,2$. The corresponding
distributions $E_1,E_2$ of $U_1\subset M$ are therefore involutive.
The involutivity of $E$ follows from Corollary \ref{gradCor}.
 Next we prove that $E_i\oplus L$
is involutive.   The formulas (\ref{e4}) show that $\n_{E_i}\tilde{X}=E_i$.
Now we check that $\n_{Y_1}Y_1'\in \Gamma (E_1\oplus L)$
for all $Y_1,Y_1'\in \Gamma (E_1)$.  Calculating the
scalar product with $Y_2\in \Gamma (E_2)$ we get
\[ g(\n_{Y_1}Y_1',Y_2)=-g(Y_1',\n_{Y_1}Y_2) = -g(Y_1',
\widehat{\n}_{Y_1}Y_2)=0,\]
since the distribution $V_2$ is parallel.\\
(ii) The fact that $E_i\oplus L\subset TM$ is parallel implies that
$E_i\subset TM_x$ is parallel.
\eprf

\noindent
Now we return to the Examples \ref{ex2} and \ref{ex3}.

In Example \ref{ex2} we have
$$V_1=T\widehat{M_1},\quad V_2=T\widehat{M_2},\quad E_1=T M_1,\quad E_2=T M_2,$$
$$X_1=\cos(s)\p_{r_1},\quad X_2=\sin(s)\p_{r_2 },\quad \a=\cos^2(s),\quad X=-\frac{1}{r}\sin(s)\cos(s)\p_s.$$
Note that $0<\a<1$.

In Example \ref{ex3} we have
$$V_1=T\widehat{M_1},\quad V_2=T\widehat{M_2}^-,\quad E_1=T M_1,\quad E_2=T M_2,$$
$$X_1=\ch(s)\p_{r_1},\quad X_2=\sh(s)\p_{r_2 },\quad \a=\ch^2(s),\quad X=-\frac{1}{r}\sh(s)\ch(s)\p_s.$$
Note that $\a>1$.

\vskip0.2cm

Let $x\in U_1$. we have two cases:
(1.) $0< \alpha(x)<1$;
(2.) $\alpha(x)<0$ or $\a(x)>1$.

{\it Case (1.)} Suppose that $0<\a(x)<1$. Then $0<\a<1$ on some open subset
$W\subset U_1$ containing the point $x$. Thus $g(\tilde{X},\tilde{X})\hat{g}(X,X)=\a -\a^2>0$ on $W$. Recall that $\tilde{X}$
is a gradient vector field, see Corollary \ref{gradCor}.
Hence we can assume that  $W$ has the form  $(a,b)\times N$, where $(a,b)\subset\bR$ and $N$ is
the level set of the function $\a$. Note also that the level sets of the function $\a$ are
integral submanifolds of the involutive distribution $E$.
Since $\t X$ is orthogonal to $E$ and  $Z(g(\t X,\t X))=0$ for all $Z\in E$, the metric $g|_W$
can be written as
$$g|_W=ds^2+g_N,$$ where $g_N$ is a family of pseudo-Riemannian metrics on $N$ depending on the parameter $s$.
 We can assume  that $\p_s=-\frac{\t X}{\sqrt{g(\t X,\t X)}}$.

By Lemma \ref{lemA6} and the Wu theorem,
the manifold $W$ is locally a product of two pseudo-Riemannian manifolds.
For $Y_1,Z_1\in E_1$ and $Y_2,Z_2\in E_2$ in virtue of Lemma \ref{lemA4}  we have
$$(L_{\t X}g)(Y_1,Z_1)=2(1-\a)g(Y_1,Z_1),\quad (L_{\t X}g)(Y_1,Y_2)=0,\quad (L_{\t X}g)(Y_2,Z_2)=-2\a g(Y_2,Z_2).$$
This means that the  one-parameter group of local diffeomorphisms of $W$ generated by the vector field $\t X$
preserves the Wu decomposition of the manifolds $W$. Hence  the manifold $(N, g\vert_N)$
can be locally decomposed into a direct product of two manifolds $N_1$ and $N_2$ which are integral
manifolds of the distributions $E_1$ and $E_2$ such that
$$g_N=h_1+h_2,$$
  where $h_i,\, i=1,2$ is a metric on $N_i$  which depends  on $s$.

From Lemmas \ref{lemA4},\ref{lemA5} it follows that the function $\a$ depends only on $s$ and satisfies the
following differential equation
$$\p_s\a=-2\sqrt{\a-\a^2}.$$
Hence, $$\a=\cos^2(s+c_1),$$ where $c_1$ is a constant.
We can assume that $c_1=0$. Since on $W$ we have $0<\a<1$ and $\p_s\a<0$, we see that
$(a,b)\subset\left(0,\frac{\pi}{2}\right)$.

Let $Y_1,Z_1\in E_1$ be vector fields on $W$ such that $[Y_1,\p_s]=[Z_1,\p_s]=0$. From (\ref{e4}) it follows that $\n_{Y_1}\p_s=-\frac{\sqrt{\a-\a^2}}{\a}Y_1$.
The Koszul formula implies that  $2g(\n_{Y_1}\p_s,Z_1)=\p_sg(Y_1,Z_1)$.
Thus we have $$-2\tg(s)g(Y_1,Z_1)=\p_sg(Y_1,Z_1).$$
This means that $$h_1=\cos^2(s)g_1,$$
where $g_1$ does not depend on $s$. Similarly,
$$h_2=\sin^2(s)g_2,$$
where $g_2$ does not depend on $s$.

For the cone over $W$ we get
$$\bR^+\times W=\bR^+\times (a,b)\times N_1\times N_2$$
and $$\widehat{g}|_{\bR^+\times W}=dr^2+r^2(ds^2+\cos^2(s)g_1+\sin^2(s)g_2).$$

Consider the functions $t_1=r\cos(s)$, $t_2=r\sin(s)$.
They define a diffeomorphism from $\bR^+\times (a,b)$ onto a subset $V\subset \bR^+\times\bR^+$.

Let $(r,y)\in\bR^+\times W\subset\widehat{M}$, then there exist a subset
$(a_1,b_1)\times (a_2,b_2)\subset V,$ where $(a_1,b_1),(a_2,b_2)\subset\bR^+$ and $r\in(a_1,b_1)$.

On the subset $$((a_1,b_1)\times N_1)\times ((a_2,b_2)\times N_2)\subset\bR^+\times W$$
the metric $\widehat{g}$ has the form
$$(dt_1^2+t_1^2g_1)+(dt_2^2+t_2^2g_2).$$

 {\it Case (2.)} Suppose that $\a(x)>1$.
Now $\p_s=-\frac{\t X}{\sqrt{\a^2 -\a }}$ and the function $\a$ satisfies
$$\p_s\a=2\sqrt{\a^2-\a}.$$
Hence, $$\a=\ch^2(s+c_1).$$
Again we can assume $c_1=0$ and from $\p_s\a>0$ we get  $(a,b)\subset\bR^+$.
For the metric $g|_{W}$ on $W=(a,b)\times N_1\times N_2$ we have
$$g|_W=-ds^2+g_N=-ds^2+ \ch^2(s) g_1+\sh^2(s)g_2.$$

The case $\a(x)<0$ is equivalent to the case $\a (x)>1$ by interchanging
the roles of $V_1$ and $V_2$, which interchanges $\a$ with $1-\a$ and $X$
with $-X$. Theorem \ref{tmain5} is proved.
\eprf

\section{Geodesics of cones}\label{geodesics}

Using Proposition \ref{doubleProp},  we now calculate the geodesics on the cone.
Suppose $\Gamma(t)=(r(t),\gamma(t))$ is a geodesic on the cone $(\widehat{M},\widehat{g})$,
where $\gamma(t)$ is a curve on $(M,g)$. Suppose we have the initial conditions
\[\Gamma(0)=(r,x)\, \text{ and }\, \dot{\Gamma}(0)=(\rho, v), \,
\]
for some $ x\in M,\, v\in T_xM $.

 Then $r(t)$ and $\gamma(t)$ satisfy
\beqa
0&=& \ddot{r}(t) -r(t) g\left(\dot{\gamma}(t), \dot{\gamma}(t)\right), \label{scalar}\\
0&=&2\ \dot{r}(t)\dot{\gamma}(t) + r(t)\n_{\dot{\gamma}(t)} \dot{\gamma}(t)\label{vector}.
\eeqa
Now one makes the following ansatz. Suppose that $\gamma$ is given as a
reparametrisation of a geodesic $\beta : \bR \rightarrow M$ of $g$:
\be\gamma(t)& =& \beta(f(t))\ee
where $\beta$ is a geodesic of $g$ with initial condition
\[\beta(0)=x\, \text{ and }\, \dot{\beta}(0)=v \neq 0,\]
implying the initial conditions for $f$:
\[f(0)=0 \, \text{ and }\, \dot{f}(0)=1.\]
As $\dot{\gamma}(t)= \dot{f}(t)\cdot \dot{\beta}(f(t))$, $g\left(\dot{\gamma}(t),
\dot{\gamma}(t)\right)= \dot{f}(t)^2g\left(\dot{\beta}(f(t)),
  \dot{\beta}(f(t))\right)$ and
$\n_{\dot{\gamma}(t)} \dot{\gamma}(t) =\ddot{f}(t)\dot{\beta}(f(t))$, we get from (\ref{scalar}) and (\ref{vector})
\beqa
0&=& \ddot{r}(t) -r(t) \dot{f}(t)^2g\left(v,v\right), \label{scalar1}\\
0&=&2\ \dot{r}(t)\dot{f}(t) + r(t)\ddot{f}(t)\label{vector1}
\eeqa
with initial conditions
\be
r(0)=r\ , && f(0)=0\ ,\\
\dot{r}(0)=\rho\ ,&& \dot{f}(0)=1.
\ee
The solution to these equations is straightforward by distinguishing several cases.

From now on we assume that $\rho\not=0$ and $v\not=0$ and consider the remaining cases for $v$ being light-like,
space-like, or  time-like.

{\em 1.) $v$ is light-like, $g(v,v)=0$, i.e. $\beta$ is a light-like geodesic.} Then the equations become
\be
0&=& \ddot{r}(t),\\
0&=&2\ \rho \dot{f}(t) + (\rho t+ r)\ddot{f}(t),
\ee
i.e.\ $r(t)=\rho t+ r$ on the one hand, and $f(t) = \frac{r t}{\rho t+r}$ on the other. This implies that $f$ and thus $\Gamma$ is defined for $t\in[0,-\frac{r}{\rho})$ if $\rho<0$, and for $t\ge 0$ otherwise.

{\em 2.) $v$ is not light-like, $g(v,v)\not=0$, i.e. $\beta$ is a space-like or time-like geodesic. } Then  we set $g(\dot{\beta}(t), \dot{\beta}(t))=g(v,v)=:\pm L^2$ with $L>0$. The equations (\ref{scalar1}) and (\ref{vector1}) become
\be
0&=& \ddot{r}(t) \mp r(t) \dot{f}(t)^2L^2, \\
0&=&2\ \dot{r}(t)\dot{f}(t) + r(t)\ddot{f}(t).
\ee
The solutions of these equations are the following
\be
r_{\pm}(t)&=& \sqrt{(\rho t+r)^2\pm L^2r^2t^2},\\
f_\pm (t)&=&\frac{1}{L} \arctan^\pm\left(\frac{L r t}{\rho t +r}\right),
\ee
in which we have introduced the notation $\arctan^+:=\arctan$ and $\arctan^-:=\mathrm{artanh}$.
Obviously $r_+$ is defined for all $t\in \rr$ whereas $f_+$ is defined for $t\in[0,-\frac{r}{\rho})$ if $\rho<0$, and for $t\ge 0$ otherwise.
The functions $r_ -$ and $f_-$ are defined on an interval $[0,T)$, where
$T$ is the first positive zero of the polynomial $((Lr-\rho)t-r)((L +r\rho)t + r)$
or $T=\infty$ if the polynomial has no positive zero.
More explicitly, $T=\frac{r}{Lr-\rho}$ if  $\rho  < Lr$ and $T=\infty$  if  $Lr \le \rho$.

\section{Cones over compact complete manifolds}
Here we generalise the proof of Gallot \cite{gallot79} for metric cones
over compact and geodesically complete pseudo-Riemannian
manifolds. We obtain the following result.

\btheo \label{tmain3} \formtmainthree \etheo
\bprf
In the Riemannian case, the values of the function $\a$
defined in \ref{e2}
are trivially restricted to the interval $[0,1]$,
since $\a = \hat{g}(X_1,X_1)\ge 0$
and $1-\a =  \hat{g}(X_2,X_2)\ge 0$.
We shall now establish the same
result for \emph{compact} complete pseudo-Riemannian manifolds $(M,g)$.
Example \ref{ex1} shows that completeness does not suffice.

\blem \label{lemA5A} Under the assumptions of Theorem \ref{tmain3},
the function $\a$ on  $M$ satisfies $0\leq\alpha\leq 1$.\elem
\bprf
On the open dense subset $U_1\subset M$ we define the vector field $$\hat X=\frac{\t X}{\sqrt{|g(\t X,\t X)|}}=\frac{\t X}{\sqrt{|\a^2-\a|}}.$$
From (\ref{e4}) and Lemma \ref{lemA5} it follows that $\n _{\h X}\h X=0$, i.e.  $\h X$ is a geodesic vector field.

Let $x\in M$ and suppose that $\a(x)<0$.
Denote by $U_x\subset U_1$ the connected component of the set $U_1$ containing the point $x$.
Since $M$ is complete, we have a geodesic $\gamma(s)$ such that $\g(0)=x$ and $\dot\g(s)=\h X(\g(s))$ if $\g(s)\in U_x$.
From Lemma \ref{lemA5} it follows that
$$\dot\g(s)\a=\h X(\g(s))\a=-2\sqrt{\a^2-\a}$$
for all $s$ such that $\g(s)\in U_x$.
Hence along the curve  $\{\g(s)|\g(s)\in \bar U_x\}$ we have
$$\a(s)=\frac{(e^{-2s}+c_1)^2}{4e^{-2s}c_1},$$
where $c_1$ is a constant. Since $\a(\g(0))<0$, we see that $c_1<0$.
If $c_1\le -1$, then for all $s>0$ we have $\g(s)\in U_x$ and $\a(\g(s))$ tends to $-\infty$ as $s$ tends to $+\infty$.
If $c_1>-1$, then for all $s<0$ we have $\g(s)\in U_x$ and $\a(\g(s))$ tends to $-\infty$ as $s$ tends to $-\infty$.
Since $M$ is compact, we get a contradiction. The case $\a(x)>1$ is similar. \eprf

Now we can  prove the  theorem  completely analogously to Gallot by verifying the same lemmas as in his proof.
\blem\label{remarque1.3}
Let $\Gamma(t)=\left(r(t), \gamma(t)\right)$ be a geodesic in $(\widehat{M},\widehat{g})$.
Then the vector field along $\Gamma$ defined by
\be H(t)&:=& r(\G(t)) \p_r(\Gamma(t)) - t\dot\Gamma(t)\ee
is parallel along $\Gamma(t)$.
\elem
\bprf The lemma follows directly from (\ref{lem1}):
\be
\widehat{\n}_{\dot\Gamma(t)}H(t)&=&
\widehat{\n}_{\dot\Gamma(t)}\left(r(t) \p_r -t \dot\Gamma (t)\right)\\
&=&\dot{r}(t)\p_r + r(t)\widehat{\n}_{\dot\Gamma(t)}\p_r - \dot\Gamma(t) - t\underbrace{\widehat{\n}_{\dot\Gamma(t)} \dot\Gamma(t)}_{=0}\\
&=&
\dot{r}(t)\p_r + r(t)\underbrace{\widehat{\n}_{\dot\gamma(t)}\p_r}_{=\frac{1}{r(t)}\dot\gamma(t)} - \underbrace{\dot\Gamma(t)}_{=\dot{r}(t)\p_r+\dot\gamma(t)}
=0,
\ee
where $r(t):=r(\G (t))$.
\eprf

\noindent
For each point $q\in \hut{M}$ we denote by $\widehat{M}^1_q$ and $\widehat{M}^2_q$  the integral manifolds of the
distributions $V_1$ and $V_2$  passing through the point $q$.
For $i=1,2$  we define the following subsets of $\widehat{M}$:
\be
C_i&:=& \left\{ p\in \widehat{M}\mid \p_r(p)\in  V_i\right\}.
\ee
Then we can prove the following lemma.
\blem\label{lemme3.2}
Let $p_1\in C_1$ (respectively, $p_2\in C_2$). Then  $\widehat{M}^2_{p_1}$
(respectively, $\widehat{M}^1_{p_2}$) is totally geodesic and flat.
\elem
\bprf
The leaves of the foliations induced by $V_1$ and $V_2$ are
totally geodesic, since both distributions are parallel.
It suffices to show that
$\widehat{M}^1_{p_2}$ is flat.

Consider a geodesic $\Gamma$ of $\widehat{M}^1_{p_2}$ starting at  $p_2=(r,x)$.
Then the vector field along this geodesic $H(t)$ defined as in Lemma \ref{remarque1.3} is parallel.
We have $H(0)=r \p_r \in V_2$ and $\dot\Gamma(0)\in V_1$ which implies
\be H(t)\  = \ r(t)\p_r - t\ \dot\Gamma(t) &\in &(V_2)_{\Gamma(t)}\text{ and }
\\
\dot\Gamma(t)&\in& (V_1)_{\Gamma(t)},
\ee
as the distributions $V_i$ are invariant under parallel transport.
Since $\p_r \inter  \hut{R}=0$, we have
\be
\hut{R}(.,.)H(t)&=& - t\ \hut{R}(.,.)\dot\Gamma(t).
\ee
Since $\hut{R}(.,.)$ are elements of the holonomy algebra leaving $V_1$ and $V_2$ invariant this implies that
\be
\hut{R}(.,.)\dot\Gamma(t)&=&0.
\ee
From this we see that the Jacobi fields along $\Gamma$ are those of a flat
manifold, which implies that $\widehat{M}^1_{p_2}$ is flat.
\eprf

Recall that we have a dense open subset $U=\{x\in\widehat{M}|\a(x)\neq 0,1 \}\subset\widehat{M}.$
\blem\label{lem4.3} Any point $p\in U$ has a flat neighbourhood. \elem

\bprf Fix a point $p\in U$.
Note that for $i=1,2$ we have $C_i\cap U=\emptyset.$

Consider the geodesic $\Gamma(t)$  starting at $p$ and satisfying the initial
condition $\dot\Gamma(0)=-r(p)X_1(p)$.
Let $H(t)$ be the vector field  along $\Gamma$ as in Lemma \ref{remarque1.3}.
We claim that if the geodesic $\Gamma(t)$ exists for $t=1$, then $\G(1)\in C_2$.
Indeed, suppose that $\Gamma(t)$ exists for $t=1$. Denote by $\tau:T_{p}\widehat{M}\to T_{\G(1)}\widehat{M}$ the parallel
displacement along $\Gamma(t)$. Since $H(1)=r(\G(1))\p_r(\G(1)) -\dot\G(1)$, we have
$r(\G(1))\p_r(\G(1))=H(1)+\dot\G(1)$. From Lemma \ref{remarque1.3}  and the fact that $\G(t)$ is a geodesic it follows that
$$r(\G(1))\p_r(\G(1))=\tau(H(0))+\tau(\dot\G(0))=\tau(r(p)\p_r(p)-r(p)X_1(p))=r(p)\tau(X_2(p))\in V_2(\G(1)).$$
This shows that $\G(1)\in C_2$.

Now we prove that the geodesic $\Gamma(t)$ exists for $t=1$. We can
apply the results of the previous section. In the notations of the previous section we have
$v=-r(p)X(p)$ and $\rho=-r(p)\a(p)$.
Since $0< \a(p)< 1$, we have  $0<L^2=g(v,v)=\a(p)-\a^2(p)$
and $r-|\rho|>0$.
Then the function $r(t)$ defining the geodesic $\Gamma(t)$  is defined on $\rr$.
The other defining function $f(t)$ is given by
\be
f(t)=\frac{1}{L}\arctan \frac{Lr(p)t}{\rho t+r(p)}=\frac{1}{L}\arctan \frac{Lt}{1-\a(p) t}.
\ee
We see that $f$ is defined for $t\in[0,1]$ as $\a(p)<1$.  Thus the  geodesic $\Gamma(t)$ is defined for $t\in[0,1]$.

Since the integral manifolds of the distribution $V_1$ are totally geodesic and $\dot\Gamma(0)\in V_1(p)$, we have
$\widehat{M}^1_p=\widehat{M}^1_{\G(1)}$.
From Lemma \ref{lemme3.2} it follows that $\widehat{M}^1_p$ is flat. Similarly we show that $\widehat{M}^2_p$
is also flat.
Hence, any  point $p$ has a flat neighbourhood.
\eprf

From  Lemma \ref{lem4.3} it follows that the dense subset $U\subset\widehat{M}$ is flat. Thus $\widehat{M}$ is flat
and $(M,g)$ has constant sectional curvature 1. This finishes the proof of
Theorem \ref{tmain3}. \eprf

\bc
If $(M,g)$ is a simply connected  compact and complete indefinite
pseudo-Riemannian manifold.  Then the holonomy algebra of the cone
$(\widehat{M},\widehat{g})$ is indecomposable.
\ec

\bprf
This follows from the fact that simply connected indefinite
pseudo-Riemannian manifolds of constant curvature are never compact.
\eprf

\section{Cones over complete manifolds}

\btheo \label{tmain4} \formtmainfore \etheo

\bprf
By going over to the universal covering, if necessary, we can assume that
$(M,g)$ is simply connected. Then $\widehat{M}$ is simply connected
and decomposable.  Let $\cup_{i\in I}W_i=U_1$ be the representation of the
open subset $U_1\subset M$ as the union of disjoint connected
open subsets. For each $W_i$ we have two possibilities: (1.) $0<\a<1$ on $W_i$; (2.) $\a<0$ or $\a>1$ on $W_i$.
Consider these two cases.

(1.) Suppose that $0<\a<1$ on $W_i$. Similarly to the proof of Theorem \ref{tmain3} we can show that the cone over $W_i$ is flat.

(2.) Suppose that $\a>1$ on $W_i$. As in the proof of Lemma \ref{lemA5A} we can show that  $\a(W_i)=(1,+\infty)$.
To proceed we need the following statement which is a generalisation of an argument used in the proof of  Theorem 27 in \cite{Bohle}.
\bs\label{PropBohle} Let $(M,g)$ be a  connected pseudo-Riemannian manifold and $\alpha\in C^\infty(M)$ with gradient $Z$ such that $g(Z,Z)\not=0$ and
$Z= (f\circ \alpha) \cdot X$ for a vector field   $X$ such that $g(X,X)=h\circ \alpha$, where
$f$ and $h$ are smooth functions on the open interval $\mathrm{Im}\  \alpha $.
If the flow $\phi$ of $X$ satisifies the following condition,
\bi
\item[($*$)] For all $c\in \mathrm{Im}\ \alpha$ exists an open interval $I_c$ such that for all $p\in F_c:=\alpha^{-1}(c)$ the interval $I_c$ is the maximal intervall on which  the flow $t\mapsto \phi_t(p)$ is defined,
\ei
then $M$ is diffeomorphic to the
product of the image of the function $\alpha$ and a level set of $\alpha$. In particular, if the manifold $(M,g)$ is geodesically complete and the vector field $X$ is a geodesic vector field, then $M$ is diffeomorphic to $\mathrm{Im}(\alpha)\times \mathrm{level\ set}$.  \es
\bprf
First we notice that
 $0\not=g(Z,Z)= (f^2\cdot h) \circ \alpha$. As $M$ is connected, we may assume that $g(Z,Z)>0$ and thus $h>0$. As the sign of $f$ plays no role in what follows we also assume that $f>0$. Furthermore, $\alpha$ satisfies the following differential equation on $M$,
 \begin{equation}\label{alphapde}
 X(\alpha)= \frac{1}{f\circ \alpha} Z(\alpha)= (f\circ \alpha) g(X,X) = ( f\cdot h)\circ \alpha.
 \end{equation}

Let $
\phi :  F_c\times I_c\rightarrow M$,
$(p,t)\mapsto \phi_t(p)$
 be the flow of the vector field $X$.
The proof is now based on the observation that if $p$ and $q$ are in the same level set of $\alpha$, then
\begin{equation}\label{observe}
\alpha (\phi_t(p)) = \alpha (\phi_s(q))\ \iff\  t=s\end{equation}
for all $t,s\in \rr$. To verify this, for each $p\in F_c$ we consider the real function
\be
\vf_c: I_c\ni t &\mapsto & \alpha(\phi_t(p))\in \mathrm{Im}\ \alpha
\ee
which satisfies the ordinary differential equation
\begin{equation}\label{ode}
\vf_c^\prime (t) = d\alpha_{\phi_t(p)}(X(\phi_t(p))= f(\vf_c (t)) \cdot h (\vf_c (t))>0.\end{equation}
Hence, for each $p\in F_c$ the function $\vf_c(t)= \alpha(\phi_t(p))$ 
is subject to the ordinary differential equation (\ref{ode}) with the same initial condition $\vf_c(0)= \alpha\left( \phi_0(p)\right)=\alpha\left( \phi_0(q)\right)=c$. Uniqueness of the solution implies that $\alpha\left(\phi_t(p)\right)=\alpha\left((\phi_t(q)\right)$ for all $t$ and all $q\in F_c$. This proves $(\Longleftarrow)$ of (\ref{observe}), and shows that $\vf_c$ does not depend on the starting point $p\in F_c$. Having this, (\ref{ode}) also  shows that 
 $\vf_c$ is strictly monotone, and thus injective which gives $(\Longrightarrow)$ of (\ref{observe}).

(\ref{observe}) shows 
that the flow $\phi$ of $X$ sends one level set $F_c$ of  $\alpha$ to another one $F_d$, i.e. $\alpha(\phi_t(p))=\alpha (\phi_t(q))$ for all $t\in I_c$ and $p,q$  in the same level set $F_c$. 

Next, we show that two level sets that are joint by an integral curve of $X$ are diffeomorphic. In fact, if $p\in F_c$ and $q=\phi_t (p)$, $\phi_t$ is a local diffeomorphism between $F_c$ and $F_d$.  $(\Longrightarrow)$ of (\ref{observe}) implies that $\phi_t|_{f_c}$ is injective. To verify that it is surjective we notice that $\phi_{-t}|_{F_d}$ is also an injective local diffeomorphism. Hence, $\phi_{-t}\circ \phi_{t}=id_{F_c}$, which implies that $\phi_t:F_c\rightarrow F_d$ is a global diffeomorphism. 

Finally, we show that for two level sets  there is at least one flow line connecting them. 
To this end, we set $\phi (F_c):=\{ \phi_t(p)\mid p\in F_c,\ t\in I_c\}$ and write
 \[ M= \bigcup_{c\in \mathrm{Im}\ \alpha} \phi(F_c).\]
We have seen that, if $F_c$ and $F_d$ are connected by an integral curve, then they are diffeomorphic under $\phi_t$.  But the maximality of  $I_c$ and $I_d$ implies that $\phi(F_c)=\phi(F_d)$. If, on the other hand, $F_c$ and $F_d$ are not joined by an integral curve then,  by maximality of $I_c$ and $I_d$, a common  point  of $\phi(F_c)$ and $\phi(F_d)$ would lie on an integral curve joining $F_c$ and $F_d$, i.e. 
$\phi(F_c)\cap\phi(F_d)=\emptyset $. In the latter case  $M$ can be written as disjoint union of open sets $\phi(F_c)$ which is not possible as $M$ was supposed to be connected.

Hence, each integral curve meets each level set once, they are all diffeomorphic, i.e. $M$ is diffeomorphic to $I_c\times F_c$. But this implies that $M$ is diffeomorphic to $(\mathrm{Im}\ \alpha)\times F$ where $F$ is a level set of $\alpha$.\eprf

Resuming the proof of the theorem we notice that the vector field $\hut{X}=\frac{\tilde{X}}{\sqrt{\alpha^2-\alpha}}$ is geodesic and proportional to the
gradient of $\alpha$. Since $(M,g)$ is complete and $\hut{X}$ is geodesic its integral curves are defined for all $t$. 
As in the proof of Proposition \ref{PropBohle} one shows that level sets are mapped onto level sets under the flow of $\hut{X}$. This shows that ($*$) in Proposition \ref{PropBohle}
 is satisfied for the vector field $\hut{X}|_{W_i}\in\Gamma (TW_i)$ on the manifold $W_i$: for $F_c$ the interval $I_c$ is limited by the real number   $a$ for which $\phi_a(F_c)\subset F_1$.
Hence, we can apply  Proposition \ref{PropBohle} to the manifolds $W_i$ and the vector field $\hut{X}|_{W_i}\in\Gamma (TW_i)$. 
Combining the result with  the proof of  Case (2.) from Theorem \ref{tmain5} yields a
 decomposition $$W_i=\bR^+\times N_1\times N_2.$$ For the metric $g|_{W_i}$ we obtain that
$$g|_{W_i}=-ds^2+ \ch^2(s) g_1+\sh^2(s)g_2,$$ where $(N_1,g_1)$ and $(N_2,g_2)$ are pseudo-Riemannian manifolds.
The fact
that the cone $(\hut W_i,\widehat{g|_{W_i}})$ is isometric to an open subset of the product of a space-like cone over
$(N_1,g_1)$ and of a time-like cone over $(N_2,g_2)$ is shown by Example \ref{ex3}.

By a variation of the proof of Theorem \ref{tmain3} we will show now that the
time-like cone over the manifold
$(N_2,g_2)$ is flat and we will explain why it is not the case for the manifold $(N_1,g_1)$.

Fix a point $p\in W_i$. Consider the geodesics $\Gamma_1(t)$ and $\G_2(t)$
starting at $p$ and satisfying the initial
conditions $\dot\Gamma_1(0)=-r(p)X_1(p)$ and $\dot\Gamma_2(0)=-r(p)X_2(p)$.

Now we prove that the geodesic $\Gamma_2(t)$ exists for $t=1$ and the geodesic $\G_1(t)$ does not exist for $t=1$.
 We can apply the results of section \ref{geodesics}. For $\G_1$ we have $v_1=-r(p)X(p)$ and $\rho_1=-r(p)\a(p)$;
for $\G_2$ we have  $v_2=r(p)X(p)$ and $\rho_2=r(p)(\a(p)-1)$.


From Section \ref{geodesics} it follows that the functions $r_1(t)$ and $f_1(t)$ defining the geodesic $\Gamma_1(t)$ are defined on the interval
$\left[0,\frac{1}{\sqrt{\a^2(p)-\a(p)}+\a(p)}\right)\subset [0,1).$
The functions $r_2(t)$ and $f_2(t)$ defining the geodesic $\Gamma_2(t)$ are defined on the interval
$\left[0,\frac{1}{\sqrt{\a^2(p)-\a(p)}-\a(p)+1}\right)\supset [0,1].$

Thus the  geodesic $\Gamma_2(t)$ is defined for $t\in[0,1]$ and the  geodesic $\Gamma_1(t)$ is not defined for all $t\in[0,1]$.

As in the proof of Theorem \ref{tmain3} we get that the manifold $\widehat{M}^2_p$ is flat. This means that
the induced connection on the distribution $V_2|_{W_i}$ is flat and  the
\emph{time-like} cone over the manifold
$(N_2,g_2)$ is flat, i.e.\ $(N_2,g_2)$  has constant sectional curvature $-1$
or $\dim N_2 \le 1$.

Note that as in Example \ref{ex1} it can be $\a>1$ on $M$, then $C_2=\emptyset$ and the induced connection on $V_1$ need not be flat.

The case $\alpha|_{W_i}<0$ is similar, with the
roles of $V_1$ and $V_2$ interchanged.
\eprf


\section{Para-K\"{a}hler cones}
\label{isosec}
\subsection*{Para-K\"{a}hler cones and para-Sasakian manifolds}
Now we consider the case where the holonomy algebra $\widehat{\mathfrak{h}}$
of the space-like cone $\widehat{M}$ over $(M,g)$ is indecomposable and preserves a
decomposition
$T_p\widehat{M}=V\+W$ ($p\in \widehat{M}$) into two complementary (necessarily degenerate)
subspaces $V$ and $W$. The next lemma\footnote{communicated to us by Lionel B\'erard Bergery}
reduces
the problem to the
case $V=V^\perp$, $W=W^\perp$.

\begin{Lem}[cf.\ Thm.\ 14.4 \cite{Tom_Krantz_PhDthesis}] Let $E$ be a
pseudo-Euclidian vector
space and
$\mathfrak{h}\subset \mathfrak{so}(E)$ an indecomposable Lie subalgebra. It $E$ admits
a non-trivial  $\mathfrak{h}$-invariant decomposition $E=V\+W$ then it admits
an  $\mathfrak{h}$-invariant decomposition $E=V'\+W'$ into a sum of totally
isotropic subspaces.
\end{Lem}
By the lemma,  we can assume that $V, W$ are totally
isotropic of the same dimension, which implies that
the metric has neutral signature.  In this
section we use a similar approach as in the previous sections but
with different structures coming up.
 These structures are related to a {\em para-complex structure}, and to a
{\em para-Sasakian structure}. We recall the basic definitions given in
\cite{cortes-mayer-mohaupt-saueressig1} and \cite{cortes-lawn-schaefer05}.

\bde
\bnum
\item Let $V$ be a real finite dimensional vector space. A {\em para-complex structure} on  $V$ is an endomorphism $J \in End(V )$, such that $J^2 = Id$ and the two eigenspaces $V^\pm :=\mathsf{ker}(Id\mp J)$ of $J$ have the same dimension.  The pair $(V,J)$ is called a {\em para-complex vector space}.
\item Let $\cal V$ be a  distribution on a manifold $M$. An {\em almost para-complex structure} on $\cal V$ is a
field
   $J\in \Gamma({End}\,\cal V)$  of paracomplex structures  in $\cal V$.
    It is called {\em integrable} or { \em  paracomplex structure  on $\cal V$} if the eigen-distributions
    $\cal V^\pm:=\mathsf{ker} (Id\mp J)$
    are involutive.
\item  A  manifold M endowed with a
para-complex structure on $TM$ is called a {\em para-complex manifold}.
\enum
 \ede
Similar to the complex case, the  integrability of $J$ is
equivalent to the vanishing of the
{\em Nijenhuis tensor} $N_J$ defined by
\beqa\label{nijenh}
N_J (X,Y)& := & J\left( [JX,Y]+[X,JY]\right) - [X,Y] - [JX,JY],\quad
X,Y\in \Gamma {\cal V}.
\eeqa

\bde
 \bnum
 \item
 Let $(V,J)$ be a para-complex vector space equipped with a scalar product $g$. $(V,J,g)$ is called
 {\em para-hermitian vector space} if $J$ is an anti-isometry for $g$, i.e.
 \beqa\label{anti}
J^*g\ :=\  g(J.,J.)& =&  -g.
 \eeqa
 \item
 A {\em (almost) para-hermitian manifold} $(M, J, g)$ is an (almost) para-complex manifold $(M, J)$ endowed with a
 pseudo-Riemannian metric $g$ such that $J^*g = -g$. The two-form $\w := g(J‡, ‡) = -g(‡, J‡)$ is called the  para-K\"ahler
 form of $(M, J, g)$.
 \item
  A {\em para-K\"{a}hler manifold} $(M, J, g)$ is a para-hermitian manifold $(M, J, g)$ such that $J$ is parallel
   with respect to the Levi-Civita-connection $\nabla$ of $g$.
\enum
 \ede
 As in the complex case, the condition $\nabla J$ is equivalent to $N_J=0$ and $d\w=0$.
In contrary to the complex case, a $2n$--dimensional para-hermitian manifold has to be of neutral signature $(n,n)$.
Note that eigen-distributions $\cal V^\pm$ of $J$ are totally isotropic and auto-orthogonal,
i.e. $(\cal V^\pm)^\bot=\cal V^\pm$. For a para-K\"{a}hler manifold the condition $\nabla J=0$ means that
the $\pm 1$-eigen-distributions $\cal V^\pm$ are  parallel. We get

\bs\label{holonomy}
A pseudo-Riemannian manifold $(M,g)$ is a para-K\"{a}hler manifold if and only if the holonomy group preserves
a decomposition of the tangent space into a direct sum of two  totally isotropic
subspaces.
\es
In the following we will show  that metric cones with  para-K\"{a}hler structure are precisely  cones over
 para-Sasakian manifolds.
\bde \label{para-SasakiDef} A {\em para-Sasakian manifold} is a pseudo-Riemannian
manifold $(M,g)$ of signature $(n+1,n)$, where $n+1$  is the
number of time-like dimensions,
endowed with a time-like geodesic unit Killing vector field $T$ such
that
$\nabla T$ defines an integrable para-complex structure
$J = \nabla T|_E : E\rightarrow E$ on $E=T^\perp$. The pair $(g,T)$ is called
a {\em para-Sasakian structure}. \ede
   Note  that the  eigen-distributions $E^\pm$  of $J = \nabla T\vert_E$ are totally isotropic and $J$ is an
   anti-isometry of $g|_E$.  Indeed, using the condition that $T$
   is a Killing field
for  $X_\pm$ and $Y_\pm$ in $\Gamma(E^\pm)$, we get
\beqa \label{isotrop}
0\ =\
( \cal L_T g)(X_\pm,Y_\pm)\
=\  g(\nabla_{X_\pm}T, Y_\pm)+  g(\nabla_{Y_\pm}T, X_\pm)\
=\  2g(X_\pm, Y_\pm).
\eeqa

A para-Sasakian manifold carries several other structures. First of all it has {\em contact structure} given by the contact form $\theta:=g(T,.)$. Indeed, for $d\theta$ we get  that
\beqa\label{dtheta}
d\theta (X_+,X_-)\ =\
-g(T,[X_+,X_-])\ =\ 2g(X_+,X_-),
\eeqa
with $X_\pm \in \Gamma(E^\pm)$. Since $E^\pm$ are dual to each other, this implies that $\theta\wedge d\theta^n\not=0$, hence $\theta$ is a contact form. The Reeb vector field of this contact structure is $T$, because
\beqa\label{dthetat}
d\theta(T,X)\ =\ -g(T, [T,X])\ =\ -g(T, \nabla_TX- \nabla_XT)\ =\ 0.
\eeqa

It also admits a
{\em para-CR structure} (see for example \cite{alek-medori-tomassini05}), which is defined on  a  $(2n+1)$-dimensional
 manifold $M$ as
an $n$-dimensional subbundle $E$ of $TM$ together with  a para-complex structure $J$ on $E$.
 For a para-Sasakian manifold this para-CR structure is given by the  one-form $\theta$. From (\ref{dtheta})
  and from the assumption that $E^\pm$ are involutive  we see that the Levi-form $L_\theta\in \Gamma( S^2E)$ of
   this para-CR structure, defined by
$L_\theta(X,Y):=d\theta|_{E}(X,JY)$, is given by the metric,
\beqa
\label{levi}
L_\theta (X,Y)\ =\ d\theta (X,JY)\ =\ -2g (X,Y)\eeqa
and is thus non-degenerate. Hence for a para-Sasakian manifold, the metric $g$ can be expressed in terms of the contact form $\theta$ and its Levi form:
\beqa\label{metric}
g\ \ =\ \ -\theta^2-\einhalb L_\theta.
\eeqa
This is in analogy to
{\em strictly pseudo-convex pseudo-Hermitian structures}
(see for example \cite{baum99} and \cite{baum99a}).
Although the definition of a para-Sasakian structure seems rather weak, it entails  the following properties.
\blem
Let $(M,g,T)$ be a para-Sasakian manifold with $E=T^\bot$ and eigen-distributions $E^\pm$. Then:
\bnum
\item $E^\pm$ are auto-parallel and $N_J|_{E^\pm}=\nabla J|_{E^\pm}=0$.
\item For $X_\pm\in \G (E^\pm)$ it holds that $\nabla_{X_-}X_+ = -g(X_+,X_-) T \mod E^+$ and \\
$\nabla_{X_+}X_- = g(X_+,X_-) T \mod E^-$.
\item For $X_\pm\in \Gamma(E^\pm)$ it holds $[T,X_\pm]\subset \Gamma (E^\pm)$.
\enum
\elem
\bprf
1. Let  $X_\pm$ and $Y_\pm$ be in $E^\pm$.
(\ref{isotrop}) implies that
\[g(\nabla_{X_\pm} Y_\pm,T)\ =\ - g(X_\pm, Y_\pm)\ =\ 0,\]
which ensures that $\nabla_{X_\pm}Y_\pm\in E$. Now, $E^\pm$ are integrable, which implies on the one hand the relation for
$N_J$, and gives on the other hand, using the Koszul formula, that $g( \nabla_{X_\pm} Y_\pm, Z_\pm)=0$ for all $Z_\pm \in
E^\pm$. Hence, $E^\pm$ are auto-parallel, which yields the relation for $\nabla J$.

2.  First of all we have that
\[g(\nabla_{X_-} X_+,T)\ =\ - g(X_+, \nabla_{X_-}T)\ =\  g(X_+,X_-).\]
Next we show that $\nabla_{X_-} X_+$ is orthogonal to $E^+$. In the following
equations $g(Y_i^-,Y_j^+)=\d_{ij}$ and
the lower indices $+,-,0$ denote the corresponding component
in $E^\pm$ and $\rr T$:
 \be
2 g( \n_{Y_i^-}Y_j^+, Y_k^+)
&\stackrel{\text{Koszul}}{=}&
g\big([Y_i^-,Y_j^+], Y_k^+\big)+ g\big([Y_k^+,Y_j^+], Y_i^-\big)+ g\big([Y_k^+,Y_i^-], Y_j^+\big)
\\
&=&
g\big([Y_i^-,Y_j^+]_-, Y_k^+\big)+ g\big([Y_k^+,Y_j^+], Y_i^-\big)+ g\big([Y_k^+,Y_i^-]_-, Y_j^+\big)
\\
&\stackrel{(\ref{dtheta})}{=}&
-\frac{1}{2}g\Big(T,
\big[
 Y_k^+, [Y_i^-,Y_j^+]_-
 \big]
 +
 \big[ [Y_k^+,Y_j^+], Y_i^-
 \big]
 +
 \big[ Y_j^+, [Y_k^+,Y_i^-]_-
 \big]
 \Big)
 \\
 &=&
 -\frac{1}{2}g\Big(T,
\big[
 Y_k^+, [Y_i^-,Y_j^+]_-
 \big]
 +
 \big[Y_i^- , [Y_j^+,Y_k^+],
 \big]
 +
 \big[ Y_j^+, [Y_k^+,Y_i^-]_-
 \big]
 \Big)
 \\
 &=&
  -\frac{1}{2}g\Big(T,
\underbrace{\big[
 Y_k^+, [Y_i^-,Y_j^+]
 \big]
 +
 \big[Y_i^- , [Y_j^+,Y_k^+],
 \big]
 +
 \big[ Y_j^+, [Y_k^+,Y_i^-]
 \big]}_{=0\text{ Jacobi identity}}
 \Big)
\\ &&
+\frac{1}{2} \underbrace{ g\Big(T,
\underbrace{\big[
 Y_k^+, [Y_i^-,Y_j^+]_+
 \big]}_{\in E^+}
 +
\underbrace{ \big[ Y_j^+, [Y_k^+,Y_i^-]_+
 \big]}_{\in E^+}
 \Big)}_{=0}
\\ &&
+\frac{1}{2}  g\Big(T,
\big[
 Y_k^+, [Y_i^-,Y_j^+]_0
 \big]
 +
 \big[ Y_j^+, [Y_k^+,Y_i^-]_0
 \big]
 \Big)
 \\
 &=&
  -\frac{1}{2}g\Big(T,
\big[
 Y_k^+,\underbrace{ g(T, [Y_i^-,Y_j^+])}_{=-2\delta_{ij}} T
 \big]
 +
 \big[ Y_j^+, \underbrace{g(T,[Y_k^+,Y_i^-])}_{=2\delta_{ik}}T
 \big]
 \Big)
 \\
 &=&
g\Big(T,
\delta_{ij} \big[
 Y_k^+,T\big] -\delta_{ik}\big[ Y_j^+,T\big]\Big)
\\ &=&0
\ee
This implies that $\nabla_{X_-} X_+\in \rr T\+ E^+$, which proves the second statement.

The last point follows from the general fact:
\begin{equation}\label{fact}
\text{If $T$ is a Killing vector  field, and $\theta=g(T,.)$, then $\cal L_T\nabla\theta = 0$.}
\end{equation}
Indeed, the Killing equation for $T$ is equivalent to $\nabla\theta =
\frac{1}{2}d\theta\in \Omega^2 M$. This implies for arbitrary tangent vectors
$X$ and $Y$ using the skew symmetry of $\nabla\theta$ that
\be
0&=&
\frac{1}{2}dd\theta (T,X,Y)
\\
&=&
T\left(\nabla\theta (X,Y)\right) -  X\left(\nabla\theta (T,Y)\right) +Y\left(\nabla\theta (T,X)\right)\\
&&- \nabla\theta ([T,X],Y) + \nabla\theta ([T,Y],X) - \nabla\theta ([X,Y],T)
\\
&=&
\left(\cal L_T\nabla \theta \right) (X,Y) -  X\left(\nabla\theta (T,Y)\right) +Y\left(\nabla\theta (T,X)\right)  - \nabla\theta ([X,Y],T)
\\
&=&
\left(\cal L_T\nabla \theta \right) (X,Y)
-  X\left(\theta (\nabla_YT)\right) +Y\left(\theta (\nabla_XT)\right)  + \theta (\nabla_{[X,Y]}T)
\\
&=&
\left(\cal L_T\nabla \theta \right) (X,Y) -\theta (R(X,Y)T)\\
&=&
\left(\cal L_T\nabla \theta \right) (X,Y).
 \ee
This can easily be applied to our situation, where we have that
\be
\nabla\theta &=& g(J.,.).
\ee
For $X_\pm\in E^\pm$ and $Y\in TM$, (\ref{fact}) implies that
\be
0& =&
\cal (L_T\n \theta )(X_\pm, Y)
\\
& =& T(g(JX_\pm,Y))-g(J([T,X_\pm]),Y)-g(JX_\pm,[T,Y])
\\
&=&
\pm \underbrace{(\cal L_Tg)(X_\pm,Y)}_{=0}\pm g([T,X_\pm],Y)-g(J([T,X_\pm]),Y),
\ee
which gives $[T,X_\pm]\in E^\pm$.
\eprf
Using these properties we obtain a description of para-Sasakian manifolds
which might look more familiar.
\bs
$(M,g,T)$ is a para-Sasakian manifold if and only if $(M,g)$ is a
pseudo-Riemannian manifold of signature $(n+1,n)$ and $T$ a time-like
geodesic unit Killing vector field, such that the endomorphism
$\phi:=\nabla T\in \Gamma(End(TM))$ satisfies:
\beqa
\phi^2&=&id + g(.,T)T \label{sasa1}\\
(\nabla_U\phi)(V)&=& -g(U,V)T + g(V,T)U,\,\,\, \forall  U,V \in TM \label{sasa2}
\eeqa
\es
\bprf
First, let $(M,g)$ be a pseudo-Riemannian manifold of signature $(n+1,n)$
with a time-like geodesic unit Killing vector field $T$ satisfying
(\ref{sasa1}) and (\ref{sasa2}). The fact that $T$ is geodesic means that
$\phi T=0$ and implies that $\phi$ preserves $E:=T^\perp$.
Putting $J:=\phi |_E$, the equation (\ref{sasa1}) shows that
$J^2=id_E$, i.e.\ $\phi$ is a skew-symmetric involution and therefore
a para-complex structure. Finally (\ref{sasa2}) ensures that $J=\phi |_E$ is
integrable because $\nabla J|_{E^\pm}=0$.

For the converse statement we assume that $(M,g,T)$ is a para-Sasakian
manifold. Setting $\phi:=\nabla T$  we get $\phi^2|_E=J^2=id$ and $\phi^2(T)=0$ which gives (\ref{sasa1}). We have to check (\ref{sasa2}): For $U=V=T$ both sides   of  (\ref{sasa2}) are zero.
For $U=X\in T^\bot$ and $V=T$ the right hand side is given by $g(T,T)X=-X$, but also the left hand side which is $(\nabla_X\phi)(T)=-\phi(\nabla_XT)=-\phi^2(X)=-X$.
For $U=T$ and $V=X_\pm \in E^\pm$ the right hand side vanishes, and the left hand side as well because of  $[T,E^\pm]\subset E^\pm$:
\be(\nabla_T\phi)(X_\pm)& =&
 \nabla_T JX_\pm - J(\nabla_T X_\pm)\\
 &=& \pm [T,X_\pm]+X_\pm -J([T,X_\pm])- J^2(X_\pm)\\
 &=&\pm [T,X_\pm] -J([T,X_\pm]) = 0.
\ee
For $U$ and $V$ both in $E^\pm$ both sides vanish because of the integrability of the para-complex structure. For $U=X_+\in E^+ $ and $V=X_-\in E^-$ the right hand side of
(\ref{sasa2}) is equal to $-g(X_+,X_-)T$ and the left hand side is given by
\[
(\nabla_{X_+}\phi) X_-\ = \ -\nabla_{X_+} X_- - \phi(\nabla_{X_+}X_-)\ =\ -g(X_+,X_-)T
\]
 because of the second point of the lemma.
\eprf

Now we can formulate the main theorem of this section.

\btheo\label{tmain6}Let $(M,g)$ be a pseudo-Riemannian manifold. There
is a one-to-one correspondence between para-Sasakian structures
$(M,g,T)$ on $(M,g)$ and para-K\"ahler structures
$(\widehat{M},\widehat{g},\widehat{J})$ on the cone
$(\widehat{M},\widehat{g})$. The correspondence is given by
$T \mapsto \widehat{J} := \widehat{\nabla}T$.\etheo

\begin{proof}
First assume that $(M,g,T)$ is a para-Sasakian manifold with
para-complex structure $J=\nabla T$ on $E:=T^\bot$, which splits
into eigen-distributions $E^\pm$. The para-complex  structure on the metric cone $(\hut{M}, \hut{g})$ is defined by
\be
\hut{J}&:=&\hut{\nabla}T.
\ee
Because of the formula for the covariant derivative of the cone, $\hut{J}$ is given by
\barr{rccclcl}
\hut{J}(\partial_r)&=&\hut{\nabla}_{\partial_r} T& =&\frac{1}{r}T\\
\hut{J}(T)&=&\hut{\nabla}_T T& =&\nabla_T T-rg(T,T)\partial_r&=&r\partial_r\\
\hut{J}(X)&=&\hut{\nabla}_X T& =&\nabla_XT-rg(X,T)\partial_r&=&J(X),\quad
X\in E,
\earr
which implies that $\hut{J}$ is an almost para-complex structure, and also an almost para-hermitian structure with respect to
the cone metric $\hut{g}$. The eigen-distributions of $\hut{J}$  are given by
\be
V^\pm&=& \rr (r\partial_r \pm T)\+ E^\pm.
\ee
They are involutive  because the distributions $E^\pm$ are involutive
and
\[\left[ r\partial_r \pm T, X_\pm\right]\ =\ \pm \left[T, X_\pm\right]\in E^\pm
\]
for $X_\pm \in \Gamma (E^\pm)$. Hence, $\hut{\nabla}T$ defines a para-K\"{a}hler structure on the cone.

 Now assume that the cone $(\hut{M},\hut{g})$ over $(M,g)$ is a
para-K\"ahler manifold with para-complex structure $\widehat{J}$.
We consider the decomposition  $T M=V^+\+V^-$ into the totally isotropic
eigen-distributions of  $\widehat{J}$.
Then the radial vector field decomposes as follows,
\beqa\label{rdecomp}
\p_r&=&\underbrace{\rho \p_r + X}_{:= X_+ \in V^+} +\underbrace{(1-\rho)\p_r -X}_{:=X_-\in V^-},
\eeqa
where  $X\in \Gamma(\hut{M})$ is a global vector field tangent to $M$. This vector field  defines  a
 para-Sasakian structure. First of all, we prove
\blem
The vector field $2rX=\widehat{J}(r\partial_r)$ on $\widehat{M}$
is tangent to $M$ and $r$-independent. Its restriction to the
submanifold
$M\cong \{ 1\}\times M\subset \bR^+\times M=\widehat{M}$  defines
is a time-like geodesic unit vector field $T$ on $(M,g)$.
\elem
\bprf
As $V^+$ and $V^-$ are totally isotropic, we get for $X$ defined in (\ref{rdecomp})
\[ 0\ =\ \rho^2 + r^2 g(X,X)\ =\ (1-\rho)^2 + r^2g(X,X). \]
This implies $\rho=\einhalb$ and $g(X,X)=-\frac{1}{4r^2}$. By the  holonomy
invariance of the distributions $V^\pm$, stated in
Proposition \ref{holonomy}, we get
\[V^+\ \ni\ \widehat{\n}_{\p_r}\left(\einhalb \p_r +X\right)\ =\ \widehat{\n}_{\p_r} X\]
and similar
\[V^-\ \ni\ \widehat{\n}_{\p_r}\left(\einhalb \p_r -X\right)\ =\ -\widehat{\n}_{\p_r} X,\]
which implies $\widehat{\n}_{\p_r} X=0$.  Hence, $[\p_r,X]=-\frac{1}{r} X$, and thus
$X=\frac{1}{2r} T$  where  $T$ is  a vector field on $M$ with  $g(T, T)=-1$. It follows that $T$
is a geodesic vector field because:
\[V^\pm
\ \ni\
\widehat{\n}_{T}\left({r} \p_r \pm T\right)
\ =\
T \pm \left(\n_{T}T + {r}\p_r \right)
\ =\ \pm\underbrace{ \left( r \p_r \pm T \right)}_{=2r X_\pm \in V^\pm} \pm \n_{T}T,
\]
i.e. $\n_{T}T\in V+\cap V^^-=\{0\}$.
\eprf

\noindent
Hence, the vector fields $X_\pm$, defined in (\ref{rdecomp}) are given by
\be
X_\pm&=&\einhalb\left( \p_r \pm \frac{1}{r} T\right)
\ee
for $T$ a time-like geodesic unit vector field on $M$. We consider now the orthogonal complement of $X_\mp$ in $V^\pm$.
\blem
Let $E^\pm:=\{Y\in V^\pm\mid \hut{g}(Y,X_\mp)=0\}\subset V^\pm$ be the orthogonal complement of $X_\mp$ in $V^\pm$. Then $E^\pm$ are tangential to $M$, orthogonal to $T$, totally isotropic, and $E:=E^+\+E^-$ is the orthogonal complement of $T$ in $TM$.
\elem
\bprf
As $V^+$ is totally isotropic any $U =a\p_r + Y\in E^+$ ($Y\in TM$)
is orthogonal to $X_+$ and $X_-$, which is equivalent to $0=a\pm rg(Y,T)$.
Hence, $a=g(Y,T)=0$. The same holds for $U\in E^-$. Both are totally isotropic
with respect to $g$ as $V^\pm$ are totally isotropic with respect to
$\hut{g}$.
\eprf

\noindent
This gives the following decomposition of the tangent bundle into three
non-degenerate distributions
\be
T\hut{M}&=& \rr \cdot \p_r \+^\bot \rr\cdot T \+^\bot\left( E^+\+E^-\right),
\ee
where $E^+$ and $E^-$ are totally isotropic.
\blem \label{8.4Lemma}
The vector field $T$ satisfies $\hut{\nabla}T|_{E^\pm}=\pm   \mathsf{id}$.
%
\elem
\bprf
The holonomy invariance of $V^\pm$ implies that
\be
\hnab X_\pm&:& T\hut{M}\rightarrow V_\pm.
\ee
But the formulae for $\hnab$ imply that
\be
\hnab X_\pm|_{TM}&=& \frac{1}{2r}\left( \mathsf{id}_{TM}\pm \hnab T\right).
\ee
Applying this to $E^\pm$ gives that $\hnab T$ leaves $E^+$ and $E^-$ invariant.
Hence, $\hnab X_\pm$ is zero on $E^\mp$, and thus
 $\hnab T|_{E^\mp}=\mp \mathsf{id}$.
\eprf
 As $T$ is a vector field on $M$, its orthogonal complement  $E$ does
not depend on the radial coordinate $r$ and defines
a distributions on M.
The same holds for $E^\pm$ because  $\hut{\n}T$ is an endomorphism on $E$ which does not
depend on $r$ and is given as $\pm \mathsf{id}$ on $E^\pm$, which are also denoted by $E$ and $E^\pm$.
 Thus, $\nabla T|_{E\pm}=\hnab T|_{E^\pm}=\pm \mathsf{id}$ defines an almost para-complex structure $J$ on
$E=T^\bot\subset TM$. As its eigen-spaces
$E^\pm$ are totally isotropic, $J$ is an anti-isometry, $g(J.,J.)=-g$.
This implies that $T$ is a Killing vector field:
\blem $T$ is a Killing vector field on $M$.
\elem
\bprf
$\left(\cal L_T g\right)(U,V)\ =\
g(J U,V) +g(JV , U)\ =\ 0
$, because $J$ is an anti-isometry.\eprf

\blem
$E^+$ and $E^-$ are involutive.
\elem
\bprf
For $Y_+$ and $Z_+$ in $E^+$ by the holonomy invariance of $V^+$, it is
$[Y_+,Z_+]\in V^+$. Hence, it suffices to show that $[Y_+,Z_+] \perp X_-$.
But this is true
because $\hnab X_-|_{E^+}=0$ (see the proof of Lemma \ref{8.4Lemma}):
\[
\widehat{g}( [Y_+,Z_+], X_-)
\ =\
-\widehat{g}( Z_+,\widehat{\n}_{Y_+} X_-) + \widehat{g}( Y_+,\widehat{\n}_{Z_+} X_-)
\ =\
0.
\]
We get the same for $E^-$.
\eprf

\noindent
Summarising we get that $T$ is a geodesic, time-like unit
Killing vector field, and $\nabla T$  is an integrable
para-complex structure on $T^\bot$.
Hence,
$(M,g,T)$ is a para-Sasakian manifold.
\eprf

\bigskip

\subsection*{Examples of para-Sasakian manifolds}
Now we  construct a family of  para-Sasakian manifolds $(M,g,T)$ of  positive non constant curvature which implies
 that the
associated  cone  $\hat M$ is not flat. We will describe $(g,T)$ in terms of
coordinates.\\
Let $(M,g,T)$ be a para-Sasaki manifold. Consider a filtration of $TM$ by  integrable distributions
$E^+\subset \rr \cdot T\+ E\subset TM$. The Frobenius Theorem implies existence of  local coordinates on $M$
  adapted to this filtration and  a hypersurface which contains the leaves  of $E^+$ and is
   transversal to $T$. We choose  local coordinates on this  hypersurface  adapted to $E^+$. Since $T$ is a Killing
   vector field, its flow can be used to extend these coordinates to coordinates
   $(t,x_1, \ldots , x_n, x_{n+1}, \ldots , x_{2n})  $ on  some open subset $ U \subset M$
   such that
\[ \ddt= T|_U\text{ and } \ddxi\in \Gamma(E^+|_U).\]
Obviously, ith respect to these coordinates  the metric $g$ is given by the matrix
of the form
\[
 \left(
 \begin{array}{c|c|c}
 -1&0&u^t \\
 \hline
 0&0&H^t \\
 \hline
 u & H & G
 \end{array}\right).
 \]
 Here  $u=(u_1,\ldots, u_n)\in C^\infty (U,\rrn)$, $H$ a non-degenerate matrix of real functions on $U$ and $G$ a symmetric
 matrix of real functions on $U$.  We choose  a basis  $Y_i^-$ of of vector fields  on $E^-$ which such that
 \be
0&=& g(T,Y_i^-),\\
\delta_{ij}&= &g(\ddxi, Y_j^-), \text{ and }\\
0&=&
g(Y_i^-, Y_j^-),
\ee
First of all, these orthogonality relations imply that
\be
Y_i^- &:=&
\left(H^{ij}u_j \right) \cdot T + b_{ij} \ddxj +H^{ij}\ddxjn
\ee
where  $H^{ij}$ is the inverse matrix to $H_{ij}$
and \beqa\label{symb}
b_{ij}+b_{ji}&=& - H^{ik}\left( u_ku_l + G_{kl}\right)H^{jl}.
\eeqa

As $T$ is a Killing vector field and $Y_i^-\in \Gamma(E^-)$, we get that  $[T,Y_i^-]=0$ which implies that
$H$, $u$, and $b$ do not depend on $t$.
Now we consider the condition (\ref{levi}) which  can be written as  $g|_E=-\einhalb L_\theta$
or
\be
-2\delta_{ij}&=&g([ \ddxi, Y_j^-], T)\\
&=& - \ddxi (H^{jk}u_k) +  \ddxi(H^{jk}) u_k
\\
&=&- H^{jk}\ddxi( u_k).
\ee
It  implies
\beqa
\label{hij}
H_{ij}&=& \einhalb \ddxj (u_i ).
\eeqa

Then  we evaluate the condition that $\n T $ acts as $-id$ on $E^-$. Note that the inverse matrix
 of the metric is given by
\[ \left(
 \begin{array}{c|c|c}
 -1&v^t&0 \\
 \hline
 v&F&H^{-1} \\
 \hline
 0 & (H^t)^{-1} & 0
 \end{array}\right),
 \]
where  $v=H^{-1}u$ and $F_{ij}=b_{ij}+b_{ji}$. We calculate
\be
\n_{Y_i^-}T &=& b_{ij}\ddxj + H^{ij} \n_{\ddxjn}T\\
&=&-v_i T + - H^{ij} \ddxjn \\
&&+ \left( b_{ij} - F_{ji} +\einhalb H^{ik}  H^{jl}\left( T(G_{kl}) +
\ddxkn (u_l) - \ddxln (u_k)\right)\right) \ddxj.
\ee
Hence, $\n_{Y_i^-}T=-Y_i^-$ is equivalent to
\be
2b_{ij} &=&  F_{ji} -\einhalb H^{ik}  H^{jl}\left( T(G_{kl}) +
\ddxkn (u_l) - \ddxln (u_k)\right).
\ee
which gives
\be
2(b_{ij} + b_{ji}) &=& 2  F_{ji} -H^{ik} H^{jl} T(G_{kl}).
\ee
This implies  that also $G$ does not depend on $t$, but together with (\ref{symb}) it also gives a formula for $b_{ij}$, namely
\beqa\label{bij}
b_{ij}&=&
-\einhalb   H^{ik}  H^{jl} \left( u_ku_l +G_{kl} + \einhalb \left( \ddxkn (u_l) - \ddxln (u_k)\right)\right).
\eeqa
Finally, we  evaluate the integrability of $E^-:=\mathsf{span} (Y_i^-)_{i=1}^n$. We write this condition as
\be
\n_{[Y_i^-,Y_j^-]} T&=& - [Y_i^-,Y_j^-]
\ee
and obtain after a lengthy but straightforward calculation that this is equivalent to
\begin{equation*}
0=\Lambda_{ij}\left[
 b_{iq}\left( \ddxq (b_{jp}) - H_{lr}b_{rp}\ddxq (H^{jl})\right)
 +H^{iq}\left(
  \ddxqn (b_{jp}) - H_{lr}b_{rp}\ddxqn (H^{jl})\right)
 \right],
\end{equation*}
in which $\Lambda_{ij}$ denotes the skew symmetrization  with respect to the indices $i$ and $j$.
Although we do not find the general solution of  this equation we
will construct  solutions with
$b_{ij}\equiv 0$. We make the following ansatz. We assume that
\be
\ddxin(u_j)&=&0,\text{ and set}\\
 G_{ij}&:=& - u_i\cdot u_j.
 \ee
 This implies that $b_{ij}=0$ which gives that
$Y_i^-= H^{ij} \left(u_j  \cdot T + \ddxjn\right)
$ with $H_{ij}=\einhalb \p_j u_i$, and
  ensures that $E^-=\mathsf{span} (Y_i^-)$ is the $(-1)$-eigen-space of $\nabla T$ and  integrable. In fact,
  we get for the Levi-Civita connection of this metric:
\be
\nabla_TT&=&0\\
\nabla_T\partial_i&=&\partial_i\\
\nabla_T\partial_{i+n}&=&-u_i T- \partial_{i+n}, \text{ i.e. }\nabla_TY_i^-=-Y_i^-\\
\nabla_{\p_i}\p_j&=&H^{kl}\p_i(H_{lj})\p_k\\
\nabla_{\p_{i+n}}\p_{j+n}&=&2 u_iu_j T+u_i\p_{j+n}+u_j\p_{i+n}\\
\nabla_{\p_i}\p_{j+n}&=&-H_{ij} T- u_j\p_i,
\ee
 which implies
 \be
\nabla_{\p_i}Y_{j}^-&=&\delta_{ij} T+ H_{kl}\p_i(H^{jk}) Y_l^-\\
\nabla_{\p_{i+n}}Y_{j}^-&=&u_i Y_j^-\\
\nabla_{Y_i^-}Y_j^-&=&0
\ee
 Now we check that the curvature of the metric is not constant.
  Calculating the curvature and  denoting  $Y_i^+$ by  $\p_i$, we
  get
\be
R( T, Y_i^\pm)  &:& \left\{
\begin{array}{rcl}
T&\mapsto & Y_i^\pm\\
Y^\mp_i &\mapsto & T
\end{array}\right.
\\
R( Y_i^\pm, Y_j^\pm)  &:& \begin{array}{rcl}
\ \ Y_k^\mp&\mapsto &\delta_{jk}Y_i^\pm -\delta_{ik}Y_j^\pm
\end{array}\\
R( Y_i^\pm, Y_j^\mp)  &:& \left\{
\begin{array}{rcl}
Y_k\pm&\mapsto &-\delta_{jk}Y_i^\pm -2 \delta_{ij} Y_k^\pm\\
Y_k\mp&\mapsto &\delta_{ik}Y_i^\mp +2 \delta_{ij} Y_k^\mp
\end{array}\right.
\ee
and the remaining terms being zero. The last terms show that $(M,g)$ does not have constant
sectional curvature. This can also be seen by calculating the derivatives of the curvature which
are zero apart from one term:
\be
\left(\nabla_{\p_i}R\right) (T,Y_j^-, \p_k, Y_l^-)&=&
-\Big( \underbrace{R(\p_i, Y_j^-,\p_k, Y_l^-)}_{=-\delta_{kj}\delta_{il}- 2 \delta_{ij}\delta_{kl}}
 +  \underbrace{R(T, Y_j^-,\p_k,\delta_{il} T)}_{=-\delta_{il}\delta_{kj}}\Big)
\\
&=& 2
\left( \delta_{kj}\delta_{il}+  \delta_{ij}\delta_{kl}\right).
\ee
Hence, $(M,g)$ is not  locally symmetric.
Note that the curvature $R$ of this metric is of the form
\[ R = R_1 -2\w\otimes J,\]
where $R_1$ is the curvature of a space of constant curvature, $J$ the para-complex structure and
 $w=g(.,J.)$ is the para-K\"{a}hler form.
Altogether we have proven:

\bs\label{psbsp}
Let $(t,x_1, \ldots, x_n, x_{n+1}, \ldots , x_{2n}) $ be coordinates on $\rr^{2n+1}$ and consider  the
metric $g$ given by
\be
 g&=& \left(
 \begin{array}{c|c|c}
 -1&0&u^t \\
 \hline
 0&0&H^t \\
 \hline
 u & H & G
 \end{array}\right)
 \ee
 in which
 \bi
 \item $u=(u_1,\ldots, u_n)$ is a diffeomorphism of $\rrn$, depending on $x_1,\ldots , x_n$,
 \item  $H = \einhalb\left(\ddxj(u_i)\right)_{i,j=1}^n$ is its  non-degenerate Jacobian,
  and
  \item $G$ is the symmetric matrix given by $G_{ij}= -u_iu_j$,
  \ei
  i.e.
  \begin{equation}\label{psmetric}
  g=-dt^2 + \sumi u_i dx_i dt + \einhalb \sumij \ddxj(u_i) dx_idx_{j+n} - \sumij u_iu_jdx_{i+n}dx_{j+n}.\end{equation}
  Then the manifold $  (\rr^{2n+1},g)$ is para-Sasakian, not locally symmetric, and its curvature is given by the following formulas
 \be
 R|_{T^\bot\times T^\bot\times T^\bot} &=&\left( R^1(J.,J.) -2\w\otimes J\right)|_{T^\bot\times T^\bot\times T^\bot},\text{ and}\\
 R(T,.)&=& R^1(T,.),
 \ee
where $R^1$ is the curvature tensor of a space of constant curvature $1$ in dimension $2n+1$, $J$ the para-complex structure and  $w=g(J.,.)$ is the para-K\"{a}hler form.
 In particular,  the space-like cone over $(\rr^{2n+1},g)$ is
para-K\"{a}hler and non-flat, i.e.\ its holonomy representation is
non-trivial and
decomposes into two totally isotropic invariant subspaces.
 \es
\bbem
\bnum
\item It is obvious that the Abelian  group $\rr^{n+1}$ acts isometrically on $(\rr^{2n+1},g)$ via
\[
\rr^{n+1}\ni(c,c_1,\ldots , c_n):\left(\begin{array}{c}
t\\
x_i\\
x_{i+n}
\end{array}\right)\mapsto
\left(\begin{array}{c}
t+c\\
x_i+c_i\\
x_{i+n}
\end{array}\right)
\]
As these isometries also fix the para-Sasaki vector field $T=\ddt$, they are automorphisms of the para-Sasaki structure $(g,T)$.
Hence, we can consider a lattice $\Gamma\subset \rr^{n+1}$ and compactify $(\rr^{2n+1},g)$ along these directions in order to obtain a para-Sasakian structure on
 \[ \rr^{2n+1}/\Gamma= T^{n+1}\times \rrn,\]
 where $T^{n+1}$ denotes the $(n+1)$-torus. We do not know under which conditions on the $u_i$'s there are more automorphisms, and if one can find enough in order to  compactify the manifold by this method.
 \item
The manifolds obtained in this way are curvature homogeneous.
\enum \ebem

More examples of para-K\"{a}hler cones are given in \cite{cortes-lawn-schaefer05} by conical special para-K\"{a}hler manifolds defined by a holomorphic prepotential of homogeneity 2. Further results on the holonomy of para-K\"{a}hler manifolds can be found in \cite{bb-ike97}.

\subsection*{Para-3-Sasakian manifolds and para-hyper-K\"ahler cones}
Now we study the case when the holonomy
algebra $\hath$ of the cone  $\widehat{M}$ preserves two complementary
isotropic subspaces $T^\pm$ and a skew-symmetric
complex structure $J$ such that $JT^+=T^-$.
Let us provide the definitions needed to formulate
a result analogous to Theorem \ref{tmain6} in this case.
\bde \label{para-3-SasakiDef}
 \bnum
\item Let $V$ be a real finite dimensional vector space.
A {\em para-hyper-complex structure} on $V$ is a triple
$(J_1,J_2,J_3=J_1J_2)$, where  $(J_1,J_2)$ is a pair of anticommuting
para-complex structures on $V$.
\item Let $M$ be a smooth manifold and $\cal V$ be a distribution on $M$.
An {\em almost para-hyper-complex structure} on $\cal V$ is a triple
$J_\a\in \Gamma({End}\,\cal V)$, $\a =1,2,3$, such that, for all $p\in M$,
$(J_1,J_2,J_3)_p$ is a para-hyper-complex structure on $\cal V_p$.  It is called
{\em integrable} if the $J_\a$ are integrable.
\item A {\em para-hyper-K\"ahler manifold} is a pseudo-Riemannian
manifold $(M,g)$ endowed with a parallel  para-hyper-complex structure
$(J_1,J_2,J_3)$ consisting of skew-symmetric endomorphisms $J_\a\in \Gamma
({End}\, TM)$.
\item
A {\em para-3-Sasakian manifold} is a pseudo-Riemannian
manifold $(M,g)$ of signature $(n+1,n)$
endowed with three orthogonal unit Killing vector fields
$(T_1,T_2,T_3)$ such that
\begin{itemize}
\item[(i)]
the vector fields $T_1$, $T_2$ are time-like and define
para-Sasakian structures $(g,T_1)$, $(g,T_2)$, see Definition
\ref{para-SasakiDef},
\item[(ii)] $T_3$ is space-like and defines
a (pseudo-)Sasakian structure $(g,T_3)$,
\item[(iii)] $\n_{T_2}T_1=T_3$,
\item[(iv)] the vector fields satisfy the following $\mathfrak{sl}_2(\bR )$
commutation relations:
\[ [T_1,T_2]=-2T_3,\quad [T_1,T_3]=-2T_2,\quad [T_2,T_3]=2T_1\quad \mbox{and}\]
\item[(v)] the tensors $\n T_\a$ define
a para-hyper-complex structure on
$E:={span}\{ T_1,T_2,T_3\}^\perp$. (Here we are using that
the conditions (ii-iii) imply
$\n_E T_\a \subset E$.)
\end{itemize}
\enum
 \ede

The assumption that $\hath$ preserves two complementary
isotropic subspaces $T^\pm\subset T_p\widehat{M}$ and a skew-symmetric
complex structure $J\in {End}\,T_p\widehat{M}$ with $JT^+=T^-$
can be now reformulated by saying that $\widehat{M}$ locally admits a
para-hyper-K\"ahler structure $(\widehat{g},\widehat{J}_1,
\widehat{J}_2,\widehat{J}_3 = \widehat{J}_1\widehat{J}_2)$, where $\widehat{J}_1|_{T^\pm}=\pm Id$ and
$(\widehat{J}_3)_p=J$. The corresponding geometry of the base manifold
$(M,g)$ is para-3-Sasakian:

\btheo\label{tmain63Sasaki}\formtmainsix3Sasaki\etheo

\bprf By Theorem \ref{tmain6} the para-Sasakian structures $(g,T_1)$
and $(g,T_2)$ induce two  para-K\"ahler structures
$(\widehat{g},\widehat{J}_1)$ and $(\widehat{g},\widehat{J}_2)$ on
the space-like cone $(\widehat{M},\widehat{g})$. Similarly, the
pseudo-Sasakian structure $(g,T_3)$
induces a pseudo-K\"ahler structure $(\widehat{g},\widehat{J}_3)$
on  $\widehat{M}$.  It suffices to show that $\widehat{J}_1
\widehat{J}_2=-\widehat{J}_2\widehat{J}_1=\widehat{J}_3$.
We recall that the vector fields $T_\a$
(considered as vector fields on $\widehat{M}$) are related to
$T_0:=r\partial_r$ by $T_\a = \widehat{J}_\a T_0$.
Using $\widehat{J}_\a = \widehat{\n } T_\a$, we show that
the conditions (iii-iv) in
Definition \ref{para-3-SasakiDef} imply that the structures
$\widehat{J}_\a$ preserve the four-dimensional distribution
$$H:=\mathrm{span}\{ T_i| i=0,1,2,3\}$$ and act as the standard
para-hyper-complex structure on $H$.  In fact, first it is clear that
$\widehat{J}_\a $ acts in the standard way on the plane $P_\a$
spanned by $T_0$ and $T_\alpha$. Second the relations (iii-iv) easily imply
that $\widehat{J_\a }$ preserves the plane $P_\a'=P_\a^\perp\cap H$.
Since $\widehat{J}_\a^2=\pm Id$,
the action of $\widehat{J_\a }$
on $P_\a'$ is completely determined by:
\begin{eqnarray*}
\widehat{J_1}T_2 &=&
\widehat{\n }_{T_2} T_1 =\nabla_{T_2} T_1
\overset{(iii)}{=}T_3\\
\widehat{J_2}T_1 &=& \n _{T_1} T_2 \overset{(iv)}{=}-2T_3+ \nabla_{T_2} T_1
=-T_3\\
\widehat{J_3}T_1 &=& \n _{T_1} T_3 \overset{(iv)}{=} -2T_2+\n _{T_3} T_1
= -2T_2-g(\n _{T_3} T_1, T_2)T_2\\ &=& -2T_2+g(T_3,\n _{T_2} T_1)T_2-T_2.
\end{eqnarray*}
This shows that the endomorphisms $\widehat{J}_\a$ act as the standard
para-hyper-complex structure on $H$.
Finally, the condition (v)   in
Definition \ref{para-3-SasakiDef}
shows that the  $\widehat{J}_\a$
act also as a para-hyper-complex structure on
$E=H^\perp$.
\eprf

\section{Lorentzian cones}

\btheo \label{tmain7}\formtmainseven\etheo

\bprf 1.  Suppose that $(\widehat{M},\widehat{g})$ admits a  parallel distribution of isotropic lines. Then there exists
on $\widehat{M}$ a nowhere vanishing  recurrent light-like vector field $p_1$. We have the decomposition $$p_1=\a\p_r+Z,$$
where $\a$ is a function on $\widehat{M}$ and $Z\in TM\subset T\widehat{M}$. Consider the open subset
$U=\{x\in\widehat{M}|\a(x)\neq 0\}$. We claim that the subset $U$ is dense in $\widehat{M}$. Indeed, suppose that $\a=0$
on an open subset $V\subset\widehat{M}$, then $p_1=Z$ on $V$. Let $Y\in TM$. We have
$$\widehat{\n}_Yp_1=\widehat{\n}_YZ=\n_YZ-rg(Y,Z)\p_r.$$ Since $p_1$ is recurrent, $$\widehat{\n}_Yp_1=\b(Y) p_1=\b(Y)
Z,$$ where $\b$ is a 1-form on $\widehat{M}$. Hence, $g(Y,Z)=0$ on $V$ for all $Y\in TM$. Thus, $Z=0$ and $p_1=0$ on $V$.

Let $y\in U$. We have $\widehat{R}_y(X,Y)\p_r=0$ for all $X,Y\in T\widehat{M}$. Hence, $\widehat{
R}_y(X,Y)p_1=\widehat{R}_y(X,Y)Z$. On the other hand, $\widehat{R}_y(X,Y)$ takes values in the holonomy algebra
$\hol_y$ and $\hol_y$ preserves the line $\bR p_{1y}$. Hence, $\widehat{R}_y(X,Y)p_{1y}=C(X,Y)p_{1y}=C(X,Y)(\a\p_r+Z)_y,$ where $C(X,Y)\in\bR$. Thus, $C(X,Y)=0$ and $\widehat{R} (X,Y)Z=0$
on $U$. Since $U\subset \widehat{M}$ is dense, $\widehat{R} (X,Y)Z=0$ on $\widehat{M}$ and $\widehat{R}(X,Y)p_1=0$ on $\widehat{M}$.

Let $x\in U$ and let $p_x=p_{1x}$. Consider any curve $\g(t)$, $t\in[a,b]$ such that $\g(a)=x$ and denote by
$\tau_\g:T_x\widehat{M}\to T_{\g(b)}\widehat{M}$ the parallel displacement along $\g$. For any $X,Y\in T_{\g(b)}\widehat{M}$ we have
$$\widehat{R}(X,Y)\tau_\g(p_x)=\widehat{R}(X,Y)(cp_{1\g(b)})=0,$$ where $c\in\bR$.
From this and the Ambrose-Singer theorem it follows that $\hol_x$ annihilates the vector $p_x$. Since $\widehat{M}=\bR^+\times M$
is simply connected, we get a parallel light-like vector field $p$ on $\widehat{M}$. Claim 1 of the theorem is proved.

Now suppose that we have a light-like  parallel vector field $p$ on $\widehat{M}$. Consider the decomposition
$$p=\a\p_r+Z,$$ where $\a$ is a function on $\widehat{M}$ and $Z\in TM\subset T\widehat{M}$. Note that
\beq\label{eL1}\widehat{g}(Z,Z)=-\a^2,\eeq and $Z$ is nowhere vanishing. As above  we can prove that the open subset
$U=\{x\in\widehat{M}|\a(x)\neq 0\}$ is dense in $\widehat{M}$.

\blem \label{lemL2} Let $Y\in TM\subset T\widehat{M}$.  We have
\begin{itemize}
\item[1.] $\p_r\a=0,$ $Y\a=rg(Y,Z)$.
\item[2.] $\widehat{\n}_{Y}Z=-\frac{\a}{r}Y$.
\item[3.] $ \widehat\n_{\p r} Z=\p_r Z+\frac{1}{r}Z=0$, i.e. $Z=\frac{1}{r}\t Z$, where $\t Z$ is a vector field on $M$.
\item[4.] $\t Z\a=-\a^2$.
\end{itemize}\elem
\bprf Claims 1-3 follow from the fact that $\widehat{\n} p=0$. Claim 4 follows from (\ref{eL1}) and  Claim 1 of the lemma. \eprf

From Claim 1 of Lemma \ref{lemL2} it follows that $\a$ can be considered as a function on $M$ and $\a$ is constant in the
directions orthogonal  to the vector field $\t Z$.

Let $x\in M$, $\a(x)\neq 0$ and let $\g(t)$ be the curve of the vector field $\t Z$ passing through the point $x$. From
Claim 1 of Lemma \ref{lemL2} it follows that along $\gamma(t)$ we have $$\a=\frac{1}{t+c},$$ where $c\in\bR$ is a
constant.

2. Suppose that $(M,g)$ is a  negative definite Riemannian manifold. In this case the vector field $\t Z$ is nowhere
light-like. From Lemma \ref{lemL2} it follows  that the gradient of the function $\a$ is equal to the vector field $\t Z$.
Hence each point  $x\in M$ has an open neighborhood $M_0$ diffeomorphic to the product $(a,b)\times N$, where $N$ is a
manifold diffeomorphic to the level sets of the function $\a|_{M_0}$. Note also that the level sets of the function
$\a|_{M_0}$ are orthogonal to the vector field $\t Z$. Consequently the metric $g|_{M_0}$ must have the following form
$$g=-ds^2+g_1,$$ where $g_1$ is a family depending on the parameter $s$ of Riemannian metrics on the level sets of the
function $\a|_{M_0}$, and $$\p_s=\frac{\t Z}{\a}.$$

From Lemma  \ref{lemL2} it follows that the function $\a|_{M_0}$ satisfies the following differential equation $$\p_s\a=-\a.$$
Hence, $$\a(s)=c_1e^{-s},$$ where $c_1\in \bR$ is a constant. Changing $s$, we can assume that  $c_1=\pm 1$. Both cases
are similar and we suppose that $c_1=1$. Note that $(a,b)=-\ln(\inf\a|_{M_0},\sup\a|_{M_0})$.

Let $Y_1,Y_2\in TM$ be vector fields orthogonal to $\t Z$ and such that $[Y_1,\p_s]=[Y_2,\p_s]=0$.
From Lemma \ref{lemL2} it follows that $\n_{Y_1}\p_s=-Y_1$.
From the Koszul formula it follows that  $2g(\n_{Y_1}\p_s,Y_2)=\p_sg(Y_1,Y_2)$.
Thus we have $$-2g_1(Y_1,Y_2)=\p_sg_1(Y_1,Y_2).$$
This means that $$g_1=e^{-2s}g_N,$$
where the metric  $g_N$ does not depend on $s$.

Thus we get the decompositions $$\widehat{M_0}=\bR^+\times (a,b) \times N$$ and
$$\widehat{g|_{M_0}}=dr^2+r^2(-ds^2+e^{-2s}g_N).$$

Define the manifold $$M_1=\bR^+\times \bR \times N$$ and extend the metric $\widehat{g|_{M_0}}$ to the metric $g_1$ on
$M_1$.

Consider the diffeomorphism
$$\bR^+\times\bR^+\to \bR^+\times\bR$$
given by $$(x,y)\mapsto\left(\sqrt{2xy},\ln\sqrt{\frac{2x}{y}}\right).$$
The inverse diffeomorphism has the form
$$(r,s)\mapsto\left(\frac{r}{2}e^s,re^{-s}\right).$$

We have
$$\begin{array}{l} \p_x=e^{-s}\p_r+\frac{e^{-s}}{r}\p_s,\\
\p_y=\frac{e^{s}}{2}\p_r-\frac{e^{s}}{2r}\p_s.\end{array}$$

We get the decomposition $$M_1=\bR^+\times\bR^+ \times N,$$
and the metric $g_1$ has the form
$$g_1=2dxdy+y^2g_N.$$

Obviously there exist two intervals $(a_1,b_1),(a_2,b_2)\subset\bR^+$ such that $1\in (a_2,b_2)$ and for
$M_2=(a_1,b_1)\times(a_2,b_2)\times N$ we have $M_2\subset \widehat{M_0}\subset M_1$. Let $g_2=\widehat{g}|_{M_2}$.

Applying Theorem 4.2 in \cite{leistner05a} to the Lorentzian situation it follows that $$\hol(M_1,g_1)=\hol(M_2,g_2)\cong\hol(N,g_N)\ltimes\bR^{\dim N},$$
where $\hol(N,g_N)$ is the holonomy algebra of the Riemannian manifold $(N,g_N)$. Thus,
$\hol(\widehat{M_0},\widehat{g|_{M_0}})\cong\hol(N,g_N)\ltimes\bR^{\dim N}.$

If the manifold $(M,g)$ is  complete, then the global decomposition follows from Proposition \ref{PropBohle}. From
Proposition \ref{completeProp} it follows that $(a,b)=\bR$ and that $(N,g_N)$ is complete. Claim 2 of the theorem is
proved.

3. Suppose that $(M,g)$ is a Lorentzian manifold. Consider the open subset $U_1=\{x\in M|\a(x)\neq 0\}\subset M$.
Obviously, $U_1$ is dense in $M$. Let $\cup_{i\in I}W_i=U_1$ be the representation of the open subset $U_1\subset M$ as a
union of disjoint connected open subsets. At each $x\in W_i$ we have $g_x(\t Z,\t Z)\neq 0$. Hence for each $W_i$ we can
use that arguments of the proof of Claim 2 of the theorem.

 Suppose that $(M,g)$ is complete. As in proof of Claim 2 of the theorem we can show that $U_i=W_i$ for each $i\in I$.
  From Claim 2 of Lemma
\ref{lemL2} it follows that the vector field $\frac{\t Z}{\a}$ is a geodesic vector field on $U_1$. Let $x\in U_1$ and let
$\g(s)$ be the geodesic such that $\g(0)=x$ and $\dot\g(s)=\frac{\t Z(\g(s))}{\a(\g(s))}$ if $\g(s)\in U_1$. Along the set
$\{\g(s)|\g(s)\in U_1\}$ we have $\a(\g(s))=e^{-s}$. Hence, $\g(s)$ is defined for all $s\in\bR$, $\g(\bR)\subset U_1$ and
$\a(\g(\bR))=\bR^+$, i.e. $(a,b)=\bR$. \eprf

\def\cprime{$'$}



\end{document}